\documentclass[EJP]{ejr} % replace ECP by EJP if needed.  
\usepackage{graphicx}
\usepackage[section]{placeins}
\usepackage[font=small]{caption}
\usepackage{amsthm}
\theoremstyle{plain}

\SHORTTITLE{Brownian Motion Moments} 

\TITLE{Mean and Variance of Brownian Motion with Given Final Value, Maximum and Argmax: Extended Version} % \thanks{}\
\AUTHORS{
  Kurt Riedel\footnote{Millennium Partners LLC
    \EMAIL{kurt.riedel@gmail.com}}}
 % Separate items with ;

\AMSSUBJ{60J65, 60J70, 91G60} % Edit. Separate items with ;
\SUBMITTED{Stochastic Models on Aug.~27, 2020, revised Feb.~21, 2021}
\ACCEPTED{May 31, 2021} % Edit.
\VOLUME{0}
\YEAR{2021}
\PAPERNUM{0}
\DOI{10.1214/YY-TN}

\ABSTRACT{The conditional expectation and conditional variance of
  Brownian motion, $B(t)$, is considered given B(t=1), its maximum and its argmax, $B(t|close, \max, \argmax)$,
  as well as those with less information: $B(t|close, \argmax)$,
  $B(t|$ $\argmax)$,  $B(t|\max, \argmax)$ where the close is the final value:
  $B(t=1)=c$ and $t \in [0,1]$. We compute the expectation and variance of a Brownian meander in time.
  By splicing together two Brownian meanders, the mean and variance of the constrained process are calculated.
  Analytic expressions are given for $E[B(t|close, \max, \argmax)]$, $E[B(t|\max, \argmax)]$ and
  $E[B(t|\argmax)]$ as well as the analogous expression for the conditional variance.
  Computational results displaying both the expectation and variance in time are presented.
Comparison of the simulation with theoretical values are shown when the close and argmax
are given. We show that using  the values of $\max$ and $\argmax$ as well as the final value, $c$,
significantly reduces the variance of the process compared to the variance of the Brownian bridge.}
\newcommand{\FSCL}{.45}
\newcommand{\argmax}{\mathrm{argmax}}

\newcommand{\close}{\mathrm{close}}
\newcommand{\BEQ}{\begin{equation}}
\newcommand{\EEQ}{\end{equation}}
\newcommand{\NEQ}{\end{equation}}
\newcommand{\Beq}{\begin{equation*}}
\newcommand{\Neq}{\end{equation*}}
\newcommand{\n}{\newline}
\newcommand{\Bt}{\tilde{B}}
\newcommand{\Th}{\Theta}

\newcommand{\phisg}{\phi_{\sigma^2}}

\newcommand{\erf}{\mathrm{erf}}
\newcommand{\erfc}{\mathrm{erfc}}
\newcommand{\Prf}{{\em Proof:}\, }
\def\distL{\stackrel{L}{\sim}}

\begin{document}
\date{ \today}

\section{Introduction}
\graphicspath{{},{momBMFig/},{/home/kriedel/Tex/momBMFig/}}
%\graphicspath{{/home/kriedel/Tex/momBMFig}}
%\graphicspath{{/home/kriedel/brownOut/pubAH/}}
The study of Brownian extrema dates back to L\`{e}vy \cite{Levy}. We give a brief overview of the results following \cite{Devroye}\  
with additional results taken from \cite{Shepp}, \cite{DI}, \cite{DIM} and \cite{ImhofHung}. Our research focuses on distribution
of $B(t)$ on $[0,1]$ given argmax of $B(t)$ and one or more of the maximum and the final value, $B(t=1)=c$.
For each Brownian path on  $[0,1]$, we denote the maximal value of $B(t)$ by $h$, denoting the ``high''. We denote the minimal value, the ``low'', by $\ell$. The location of the first time when the $B(t)$ reaches $h$ will be denoted by $\theta$ (and occasionally by $argmax$). The location of the first time when the $B(t)$ reaches $\ell$ will be denoted by $\theta_{\ell}$.
We are interested in the density of Brownian paths conditional on $(c,h,\ell, \theta,\theta_{\ell})$. Unfortunately, this five dimensional parameterization is too complex for us to evaluate. Instead, we settle for evaluating the conditional density, $p(x,t|c, \theta, h)$ as well as $p(x,t|\theta, h)$ and  $p(x,t|\theta)$.
The distribution, $P(x,t |c, \theta, h)$, is well known \cite{WilliamsBAMS,Denisov},
and has a representation in terms of two Brownian meanders \cite{DI,DIM, ImhofJAP,ImhofHung, Biane}.
We do give a more explicit representation of this density in Corollary \ref{denWill} .
Our companion article \cite{Riedel} examines on $p(x,t|c, h,\ell)$, $p(x,t|c, h)$, $p(x,t|h,\ell)$ and $p(x,t|h)$.

To better understand the conditional probability, we calculate the first and second moments of the probability conditional with respect to $x$. We evaluate the distribution, mean and variance of $B(t|close, \argmax, max)$, $B(t|\max, \argmax)$ and $B(t|\argmax)$. Our companion article \cite{Riedel} examines on $B(t|\max)$, $B(t|close, \max)$,  and $B(t|close, \max, \min)$.
For our results, we view the maximum value, the final value and the location of the maximum, $\theta$, as statistics on
the Brownian motion. In other words, given values of $(c,h,\theta)$, we can estimate the properties of a realization of Brownian motion
by plugging in the values of these statistics into the density and its first and second moments.
 The calculation of the moments of $p(x,t| c,\theta, h)$
is new, as is the evaluation of   $p(x,t|\theta, h)$ and  $p(x,t|\theta)$ and their moments.

By computing the moments, $E[B(t|close, \max, \argmax)]$ and $Var[B(t|close, \max, \argmax)]$, we show how
the values of (close, high, low, argmax) determine
the behavior of $B(t)$ over the entirety of the  region, $[0,1]$. Both the analytic formulae
of the conditional expectation and variance as well as the figures in this paper build intuition about Brownian paths.

Section 2 reviews the distribution of Brownian motion extrema. Section \ref{BMsect} reviews the Brownian meander and calculates its expectation and variance in Theorem \ref{BMmom}. Section \ref{SBMsect} elaborates on Williams's construction \cite{WilliamsBAMS,WilliamsProc,Denisov} of conditional Brownian motion as the splicing together of two Brownian meanders. Theorem \ref{BrownMom} gives analytic
expressions for the expectation and variance of the conditional density, $E[B(t|c,\theta,h)]$ and $Var[B(t|c, \theta,h)]$.
%Section \ref{EVCA} presents the analytic and simulation results for $B(t|close, \argmax)$.
Section \ref{EHThsect} calculates the density (Theorem \ref{denTH}) and moments, $E[B(t|\theta,h)]$ and $Var[B(t|\theta,h)]$, of $B(t|\theta=\argmax,h=\max)$.
Section \ref{SimDesc} and Appendix B describe our simulation framework and compares the simulation with our analytic expressions for the conditional mean and variance.
%We compare the analytic moments with our simulation.
Section \ref{EVAsect} presents our analytic $B(t|\argmax)$ including expressions for $E[B(t|\theta)]$ and $Var[B(t|\theta)]$.
It compares our analytic moments with those of the simulation. 
We conclude by comparing the time averaged expected variance given one or more of $({\rm close}, \argmax, \max, \min )$.
%Section \ref{EVCAH}
In this extended version, we add many more figures of the $E[B(t|{\mathrm parameters})]$ and $Var[B(t|{\mathrm parameters})]$.
Appendix C presents the simulation results for the moments of $B(t|\max, \argmax)$.
We display plots of $E[B(t|\theta,h)]$ and $Var[B(t|\theta,h)]$ as $\theta$ and $h$ vary.
Appendix D presents similar simulation results for the moments of $B(t|close, \argmax, \max)$.
Appendix E concludes simulation results for the moments of $B(t|close, \argmax)$.
%Section \ref{MeanderComp} compares the simulation results with the analytic representation as moments of two spliced meanders.
%Section \ref{EVCAHav} compares time averages of $E[B(t)]$ and $Var[B(t)]$ versus the givens.

% Appendix C displays plots of $E[B(t|\theta,h)]$ and $Var[B(t|\theta,h)]$ as $\theta$ and $h$ vary.
% Appendix D displays plots of $E[B(t|\theta,c)]$ and $Var[B(t|\theta,c)]$ as $\theta$ and $c$ vary.

\section{Distributions of Brownian Extrema}

%corresponding to the distribution $H_c \distL .5*( c + \sqrt{c^2 +2 E})$ where $E$ is a standard exponential random variable.

The study of Brownian extrema date back to the founders of the field \cite{Levy}. Our brief discussion follows \cite{Devroye} with additional results taken from \cite{Shepp,DI,DIM,ImhofHung,Borodin}. %See also \cite{BertoinPitman,KS,PY,Yor}.
%#Independent unit Gaussian variables will be denoted by $N$ and $N'$. Exponentially distributed random variables will be denoted by $E$ and $E'$.
To our knowledge, there is no known expression for $p(c,h, \ell, \theta, \theta_{\ell})$. %The distribution of b
There are expressions for $p(c,\theta,h)$, $p(c,h, \ell)$, $p(c, h)$, $p(h, \ell)$, $p(h,\theta)$, $p(c,\theta)$ as well as
$p(c)$, $p(h)$ and $p(\theta)$.

In \cite{Shepp,KS}, the joint density of the close, $B(t=1)\equiv c$, the high, $\max_{ \{ t\le 1\}}B(t)\equiv h$, the location of the high, $\argmax_{ \{ t\le 1\}}B(t)\equiv \theta$, is shown to be:
\begin{align} \label{SheppResult0}
  p(c,\theta,h)&= \frac{h(h-c)}{\pi \theta^{3/2}(1-\theta)^{3/2}}
\exp\left(-\frac{h^2}{2\theta}-\frac{(h-c)^2}{2(1-\theta)} \right) \ , \ h\ge 0, \ h \ge c.
 \end{align}
Conditional density of $\theta$ and $h$ given $B(1)$ is
\BEQ \label{SheppCond}
p(\theta, h|B(1)=c) = \frac{h(h-c)\sqrt{2\pi}\exp(c^2 /2)}{\pi \theta^{3/2}(1-\theta)^{3/2}}
\exp\left(-\frac{h^2}{2\theta}-\frac{(h-c)^2}{2(1-\theta)} \right) \ , \ h\ge 0, \ h \ge c.
\NEQ
The unconditional density of $\theta$ and $h$ is obtained by integrating \eqref{SheppResult0} on $c < h$ \cite{Shepp,Devroye,AC}:
\BEQ \label{eq:denThH}
p(\theta,h) =\frac{1}{\pi}\frac{h \exp(-h^2/2\theta)}{\theta^{3/2}(1-\theta)^{1/2}} \ , \ h>0, \ 0< \theta<1.
\NEQ
%The corresponding conditional density is $ p(h|\theta) = (h/\theta) \exp( -h^2 /2\theta)$ \cite{Shepp}. The expectation and variance satisfy:
The conditional density of  $h$ given $\theta$ \cite{Shepp} is 
\BEQ \label{eq:condenThH}
p(h|\theta) =\frac{h}{\theta}\exp(-h^2/2\theta) \ , \ h>0.
\NEQ
Thus $E[h|\theta] = \sqrt{\pi \theta/ 2}$, $Var[h|\theta] = (2-\pi/2) \theta$.
The classic result \cite{Levy,ItoMcKean} derived using the reflection principle is
the joint distribution of the close, $c$, the high, $h$ is
\BEQ \label{eq:disHC0}
P(h,c)= P(\max\{B(s),\ s \in [0,1]\}\le h, B(1)=c) = \phisg(c) - \phisg(2h-c) \ 
\NEQ
where we denote the Gaussian density by $\phi_{\sigma^2}(x) = (2\pi \sigma^2)^{-.5} \exp(-x^2/2\sigma^2)$.
This result has been generalized to other diffusions in \cite{PY}.
The marginal density of the maximum, $h$, and $B(1)$: 
\BEQ \label{eq:denHC0}
p(h,c) =\sqrt{\frac{2}{\pi}} (2h-c) \exp(-(2h-c)^2/2) \ , \ h\ge 0, \ h \ge c \ .
\NEQ
\cite{Levy,ItoMcKean}. This can be derived by differentiating \eqref{eq:denHC0} or  by integrating \eqref{SheppResult0} over $\theta$.
%\cite{WilliamsBAMS,ImhofJAP}.  %{BertoinPitman}

%corresponding to the distribution $H_c \distL .5*( c + \sqrt{c^2 +2 E})$ where $E$ is an exponential random variable.

%We denote the Gaussian density by $\phi_s(x) = (2\pi s)^{-.5} \exp(-x^2/2s)$.
%Independent unit Gaussian variables will be denoted by $N$ and $N'$. Standard exponentially distributed random variables will be denoted by $E$ and $E'$.
%The density of the high, (maximum of $B(t)$), $h$ is that of the half normal: $2|\phi_1(x)|$. 
%The distribution of the high, given the closing value, $B(1)=c$, is 
%%As discussed in Devroye, the distribution of $M$ given $B(1)=r$ is
%\BEQ \label{eq:HcondC}
%F(h|c) = 1 -\exp\left((c^2 -(2h-c)^2)/2 \right) \ .
%\NEQ
Let $H_c$ be the random variable $\max_{\{0<t\le 1\}}  B(t)$  conditional on $B(t=1)=c$.
The density of $H_c$, $p(h|c)$,  can be computed using \cite{McLeish,Devroye}:
%The density for \eqref{eq:HcondC}  can be computed using \cite{McLeish,Devroye}:  
\BEQ \label{eq:BP}
H_c \equiv \max B(t|B(1)=c) \distL c/2 +  \sqrt{c^2 +2 E} /2
\NEQ
where $E$ is a standard exponentially distributed random variable and $\distL$ means equality in distribution.
The marginal density of $h$ is the half normal: %$2x$ of the normal density:
\BEQ
p(h) =\sqrt{\frac{2}{\pi }} \exp(-h^2/2) \ ,\ h>0.
\NEQ
%The density of the high, (maximum of $B(t)$), $h$ is that of the half normal: $p(h)=2 \phisg(h) =\sqrt{\frac{2}{\pi \sgp }} \exp(-h^2/2\sgp)$, $h>0$.
The marginal density of the location, $\theta$, of the maximum of $B(t)$ is the well-known
arcsine distribution,
$p(\theta) = \frac{1}{\pi } \frac{1}{\sqrt{\theta (1-\theta)}}$ . 
%\BEQ\NEQ
Using \eqref{SheppResult0}, one can generate the joint distribution of $(\theta, M,B(1))$ \cite{Devroye} as 
\BEQ \label{eq:ThMBrep}
(\theta, M,B(1))\distL \left(\Th \equiv \frac{1+ cos(2\pi U)}{2}, \sqrt{2\Th E}, \sqrt{2\Th E}-\sqrt{2(1-\Th) E'} \right) \ \ .
\NEQ
An immediate consequence of \eqref{eq:ThMBrep} is that the variance of $B(1)$ given $\theta$, the time of the maximum of $B(t)$, is
independent of $\theta$ and satisfies:
\BEQ \label{eq:B1condTh}
Var[ B(1) | \theta] = 2(1 - \pi/4) \ ,\quad E[ B(1) | \theta] = \sqrt{2}\left(\sqrt{\theta} -\sqrt{1-\theta}\right) \ \ .
\NEQ
Here $E$, $E'$, $U$ are independent random variables where $E$ and $E'$ are standard exponentially distributed and $U$ is uniformly distributed in $[0,1]$.
%N.sqrt(2-N.pi/2.)=0.65514
We are unaware of any previous derivation of this straightforward result. We actually discovered it during our simulation.
%More generally, the unconditional densityt of $\theta$ and $c$ is obtained by integrating \eqref{SheppResult0} on $(c < h)\cap (h\ge 0)$:
%BELOW NOT DONE:
%\BEQ \label{eq:denThH}
%p(\theta,h) =\frac{1}{\pi}\frac{h \exp(-h^2/2\theta)}{\theta^{3/2}(1-\theta)^{1/2}} \ , \ h>0, \ 0< \theta<1.
%\NEQ
%where we use
%\BEQ
%\int_0^{\infty} h(h-c) \exp(-\frac{(x-\kappa)^2}{2\sigma^2})= \frac{\sigma^4}{2}\left[\sqrt{2\pi}\sigma(1+ \sigma^2\kappa^2) \left(1+\erf(\frac{\kappa}{\sqrt{2}\sigma}+\frac{2\kappa}{\sigma^2} e^{\frac{-\kappa^2}{2\sigma^2}}\right]-c\left[\sigma^2 e^{\frac{-\kappa^2}{2\sigma^2}}  \right]
%\NEQ
We now evaluate $p(\theta,c)$. For $c>0$:
\begin{align} \label{eq:denThCp}
p(\theta,c)= \int_c^{\infty} p(\theta,h,c)dh = \int_c^{\infty} \frac{h(h-c)}{\pi \theta^{3/2}(1-\theta)^{3/2}}
\exp\left(-\frac{h^2}{2\theta}-\frac{(h-c)^2}{2(1-\theta)} \right)dh 
 \end{align}
\begin{align} \label{denThCpA}
=\  \frac{c \theta  e^{\frac{-c^2}{2\theta}}}{\pi \sqrt{\theta(1-\theta)}}  -\frac{(c^2-1)e^{\frac{-c^2}{2}} \erfc(|c|\sqrt{\frac{(1-\theta)}{2\theta}})}{\sqrt{2 \pi}}  \ .
\end{align}
Here $\erfc$ is the complimentary error function, $\erfc(z) = 1 - \erf(z)$.
%integrate x(x-c) \exp(- x^2/(2*.93)) \exp(-(x-c)^2/.14) from  c to infinity
%where $\kappa_{\theta} = \theta/(1-\theta)$.
For $c\le 0$, the same integral over $h$ runs from $0$ to $\infty$.
\begin{align} \label{eq:denThCn}
p(\theta,c)= 
\frac{|c| (1-\theta) e^{\frac{-c^2}{2(1-\theta)}}}{\pi \sqrt{\theta(1-\theta)}}-\frac{(c^2-1)e^{\frac{-c^2}{2}} \erfc(|c|\sqrt{\frac{\kappa_{\theta}}{2}})}{\sqrt{2 \pi}} \ \ .
%- \sqrt{2 \pi \theta(1-\theta)}(c^2-1)e^{\frac{-c^2}{2}} \erfc(|c|\sqrt{\frac{\theta}{2(1-\theta)}})  \right]
\end{align}
This result corresponds to the conditional density $p(\theta|c)$ in \cite{AC}.

A result that goes back to L\`{e}vy \cite{Levy, ChoiRoh}, if not earlier,  is 
\begin{theorem}
The joint distribution of the close, $c$, the high, $h$, and the low, $\ell$ is
\BEQ \label{ChoiRoh0}
P(c,h,\ell)=P(B(1)=c, \ell \le \{B(s),\ s \in [0,1]\}\le h) = \sum_{k=-\infty}^{\infty} \phisg(c-2k(h-\ell) ) - \phisg(c-2h- 2k(h-\ell) 
\NEQ
\BEQ \label{ChoiRoh1}
=  \phisg(c) - \sum_{k=0}^{\infty} \left[\phisg(c-2h- 2k\Delta) + \phisg(c-2\ell+ 2k\Delta) \right]+ \sum_{k=1}^{\infty}\left[\phisg(c-2k\Delta)+ \phisg(c+2k\Delta)\right]
\NEQ
where $\Delta \equiv (h-\ell)$.
\end{theorem}
Although \eqref{ChoiRoh0} is classic,  \eqref{ChoiRoh1} seems to be new.
By refactoring \eqref{ChoiRoh0} into the symmetric form, \eqref{ChoiRoh1}, we not only treat $h$ and $\ell$ symmetrically, but also show that the series is in an alternating form. Here $P(c,h,\ell)$ is a distribution in $h,\ell$ and a density in $c$.
The density, $p(h,\ell)$ may be found in \cite{Riedel}.

\section{Brownian Meander and Its Moments} \label{BMsect}

A Brownian meander, $B^{me}(t)$, can be thought of as a Brownian motion restricted to those paths where $B^{me}(t) \ge 0$.
Since this is a set of paths with measure zero, we rescale the Brownian paths to get a Brownian meander. Following \cite{Belkin,DIM},
let $B(t)$ be a Brownian motion, $\tau_1= \sup\{t \in [0,1]:\ B(t)=0\}$, $\Delta_1=1-\tau_1$,
we define the Brownian meander as
$B^{me}(t) = |B(\tau_1 + t \Delta_1)| / \sqrt{ \Delta_1}$.
By $B^{me}_c(t)$, we denote the Brownian meander restricted to $B^{me}(t=1) =c$.
Results for the distribution of meanders can be found in \cite{DI,DIM, ImhofJAP,ImhofHung, Biane}.
%We review two representations for the meander.
Brownian meanders at a fixed time, $t$, are distributed as  \cite{ImhofJAP, Devroye}:
\BEQ 
\label{eq:BPA}
B^{me}_c(t) \distL \sqrt{\left(ct+\sqrt{t(1-t)}N\right) + 2Et(1-t)} \ \ ,
\NEQ
where $N$ is a standard normal variable, $E$ is a standard exponentially distributed random variable, and $N$ and $E$ are independent. Equation \eqref{eq:BPA} applies for a single fixed time and
not for the process, $B(t|B(t=1)=c)$ for all $t \le 1$.
%In \cite{DIM}, From DIM  Meander, $B^{me}(t)$.
%Let $B(t)$ be the standard Brownian Motion,
%Let $\tau= \sup \{t \in [0,1], W(t) = 0 \}$, $\delta_1=1-\tau$. Define $B^{me}(t)= \delta_1^{-1}|W(\tau+ t\Delta_1)|$ for $t \in [0,1]$.
We are more interested in the transition density of the meander.
Following \cite{DIM}, let $\phi_s(x) = (2\pi s)^{-.5} \exp(-x^2/2s)$, $N_x(a,b)= \int_a^b \phi_s(x)dx$, then

\begin{theorem}\cite{DIM,ImhofJAP}
$B^{me}$ has transition density: $p(B^{me}(t)=y) = 2 y t^{-3/2}\exp(-y^2/2t)N_{1-t}(0,y)$ for $0\le t \le 1$, $y>0$.
For $0<s< t \le 1$, $x,y>0$,
\BEQ \label{eq:Meander0}
  p(B^{me}(t)=y|B^{me}(s)=x) = g_{t-s}(x,y) N_{1-t}(0,y)/N_{1-s}(0,x) \ ,
\NEQ
\BEQ \label{eq:Meand2}
g_t(x,y) \equiv \phi_t(y-x) -\phi_t(y+x) \ .
%P(B^{me}(t)=y|B^{me}(s)=x) = g(t-s,x,y) N_{1-t}(0,y)/N_{1-s}(0,x)
\NEQ
$P(B^{me}(1)\le x) = 1 -\exp(-x^2/2)$ is the Rayleigh distribution.
\end{theorem}
In \eqref{eq:Meander0}, we slightly abuse notation: The probability that $B^{me}(s)=x$ is infinitesimally small and is a density in $x$: $p(B^{me}(t)=y|B^{me}(s)= dx) = p(B^{me}(t)=y|B^{me}(s)=x) dx$.
With this understanding, we will continue to write the conditional densities without $dx$.
The Bayesian Rule gives the conditional density given the final value of a Brownian meander:
\begin{corollary}\cite{Devroye}
 For $0<s< t \le 1$, $x,c>0$,  $B^{me}$ has density:  
\BEQ %\begin{align}
  \label{eq:BayesMeander}
  p(B^{me}(s)=x|B^{me}(t)=c)   =\   p(B^{me}(t)=c|B^{me}(s)=x)  p(B^{me}(s)=x)/p(B^{me}(t)=c)
\NEQ
\Beq
 =\ g_{t-s}(x,c) \frac{x t^{3/2}}{cs^{3/2}} \exp(-x^2/2s) \exp(c^2/2t) \ = g_{t-s}(x/\mu,\mu c) \frac{x t^{3/2}}{cs^{3/2}}  \  
\Neq %\end{align}
where $\mu \equiv \sqrt{s/t}$. For $t=1$, \eqref{eq:BayesMeander} can be rewritten as $p(B^{me}(s)=x|B^{me}(t)=c) = g_{s(1-s)}x,c'=sc) \frac{x}{c'} $.
\end{corollary}
Equivalent formulae  to \eqref{eq:BayesMeander} are found in \cite{Devroye}.
In Appendix A, we calculate the expectation and variance of the Brownian meander using this transition density:
\begin{theorem} \label{BMmom}
  For the  Brownian meander, $B^{me}(t)$, the first and second moments are
\begin{align} 
  M_1(s,t,c) &= \int_0^{\infty}x p(B^{me}(s)=x|B^{me}(t)=c) dx \label{mnMeander1}\\
  &= \frac{[t-s+ sc^2/t] {\rm erf}(\frac{c\sqrt{s}}{\sqrt{2t(t-s)}})}{c} +\frac{\sqrt{2s(t-s)}}{\sqrt{\pi t}} \exp(-sc^2/2t(t-s)) \ ,
\end{align}
\BEQ
  M_2(s,t,c) = \int_0^{\infty} x^2 p(B^{me}(s)=x|B^{me}(t)=c) dx 
   = \frac{3s(t-s)}{t}+ c^2s^2/t^2 \ . \label{mnMeander2}\\
   \NEQ
\end{theorem}
\Prf These formulae are derived in Appendix A using the  transistion density \eqref{eq:BayesMeander}.

Theorem \ref{BMmom} is new and is the basis for the moment calculations in the remainder of the article.
If $t=1$, we often suppress the third argument: $M_1(s,c)=M_1(s,c,t=1)$. The variance of the meander satisfies:
$Vr(s,c)=M_2(s,c,t=1)-M_1(s,c,t=1)^2$. Thus, $Vr(s,c)$ is the variance of a Brownian
meander when the end time, $t=1$ with final value $c$.

\section{Spliced Brownian Meander Representation and Its Moments} \label{SBMsect}
Following Williams, \cite{WilliamsBAMS,WilliamsProc,Denisov,Devroye}, we represent a Brownian motion on $[0,1]$ 
given the location, $\theta$, its maximum, $h$ and its final value $c$.
We represent $B(t)$ in $[0,\theta]$ in terms of a meander $B^{me}(t)$:
\begin{theorem}\label{CHAThm}
  Consider a Brownian motion  conditional on
  $h=\max_{0\le s \le 1}\{B(s)\}$, $\theta$ is the smallest $s$ such that $B(\theta)=h$ and $B(1)=c$.
For a fixed time $t$, $B(t|c,\theta,h)$ is distributed as a scaled Brownian meander to the left of $\theta$ and a second scaled Brownian meander to the right of $\theta$.
\BEQ \label{CHAleft}
B(t|c,\theta,h) \distL h - \sqrt{\theta} B^{me}_r (1 -t/\theta) \ \ {\rm for}\ \  t \le \theta \ \ {\rm where}\ \ r \equiv h / \sqrt{\theta} \ \ .
\NEQ
Here $B^{me}_r (t)$ is the meander process with $B^{me}_r(t=1) =r$ and $r \equiv h / \sqrt{\theta}$.
Similarly for $t \in [\theta,1]$, define 
\BEQ \label{CHAright}
B(t|c,\theta, h) \distL  h - \sqrt{1-\theta} \Bt^{me}_q (\frac{t-\theta}{1-\theta} ) \ \ {\rm for}\ \ t \ge \theta \ \   {\rm where}\ \ q \equiv \frac{h-c}{ \sqrt{1-\theta}}
\NEQ
and $\Bt^{me}_q$ is an independent Brownian meander.
\end{theorem}
The formal proof may be found in \cite{WilliamsBAMS,WilliamsProc,Denisov}. The result is easy to understand: A Brownian meander, $B^{me}_c(t)$, is equivalent to a Brownian bridge restricted to non-negative paths. On either side of the first maximum, $\theta$, the Brownian paths are restricted to not go above the maximum.  Thus on either side of $\theta$, we subtract off the appropriately scaled Brownian meander.
%On either side of the first maximum, we subtract of the appropriately scaled meander.

A rigorous notation would replace $ B(t|c,\theta,h)$ with $B(t|dc,d\theta,dh)$.
Theorem \ref{CHAThm} is well known and sometimes called the Williams representation.
We need a more explicit representation to calculate the moments. We now derive it
using the meander density representation of the previous situation. We believe the following corrolary is
more explicit than previous analyses.
%\begin{align} \label{BMDensLo}
%P(B(s)=x|\theta,h,c) = g_{s/\theta}(\frac{h-x}{\sqrt{\theta}},\frac{h}{\sqrt{\theta}}) \frac{\sqrt{\theta}(h-x)}{(1-\frac{s}{\theta})^{3/2}h} \exp(\frac{h^2}{2\theta}) \exp(\frac{-(h-x)^2}{2(\theta-s)})  \ = %{\rm for}\ s< \theta.
%\end{align}
\begin{corollary} \label{denWill}
  % The density of a  Brownian motion conditional on $h=\max_{0\ge s \ge 1}\{B(s)\}$, $\theta} is the smallest $s$ such that $B(\theta)=h$ and $(B(1)=c}$
Under the assumptions of Theorem \ref{CHAThm}, the conditional density of a Brownian motion at a fixed time, $t$, satisfies
\begin{align} \label{BMDensL}
p(B(t)=x|c,\theta,h) = g_{t}( h-x,h) \frac{\theta^{1.5}(h-x)}{(\theta-t)^{1.5}h} \exp(\frac{h^2}{2\theta}) \exp(\frac{-(h-x)^2}{2(\theta-t)}) \ \ {\rm for}\ \ t< \theta \ \ .
\end{align}
For $t> \theta$:
%\begin{align} \label{BMDensRo}
%  P(B(s)=x|\theta,h,c) = g_{\frac{1-s}{1-\theta}}(\frac{h-x}{\sqrt{1-\theta}},\frac{h-c}{\sqrt{1-\theta}}) \frac{(1-\theta)^2(h-x)}{(s-\theta)^{\frac{3}{2}}(h-c)} \exp(\frac{-(h-x)^2}{2(s-\theta)}) \exp(\frac{(h-c)^2}{1-\theta)})  \ .
%\end{align}
\begin{align} \label{BMDensR}
  p(B(t)=x|c,\theta,h) = g_{1-t}(h-x,h-c) \frac{(h-x)(1-\theta)^{1.5}}{(h-c)(t-\theta)^{1.5}} \exp(\frac{-(h-x)^2}{2(t-\theta)}) \exp(\frac{(h-c)^2}{2(1-\theta)})   \ \ .
\end{align}
For  $t \le \theta$, the joint density satisfies
\begin{align} \label{BMDensXThHC}
  p(B(t)=x,\theta,h,c)= \frac{(h-c)(h-x)}{\pi(1- \theta)^{1.5} (\theta-t)^{1.5}}  g_{t}(h-x,h)
\exp\left(\frac{-(h-x)^2}{2(\theta-s)} -\frac{(h-c)^2}{2(1-\theta)} \right)
\end{align}
where $t \le \theta$. For $t>\theta$, the joint density is
%we multiply \eqref{BMDensR} by $p(\theta,h,c)$: % from \eqref{SheppResult0}
\begin{align} \label{BMDensXThHCR}
  p(B(t)=x,\theta,h,c)= \frac{h(h-x)}{\pi \theta^{1.5} (t-\theta)^{1.5}}  g_{1-t}(h-x,h-c)
  \exp(\frac{-(h-x)^2}{2(t-\theta)}-\frac{h^2}{2\theta}) \ \ .
%\exp\left(-\frac{h^2}{2\theta}-\frac{(h-c)^2}{2(1-\theta)} \right)
\end{align}
\end{corollary}
\Prf To prove \eqref{BMDensL} and \eqref{BMDensR}, we insert the transition density of \eqref{eq:BayesMeander}
into the results of Theorem \ref{CHAThm}.
To prove \eqref{BMDensXThHC}, we  multiply \eqref{BMDensL} by $p(\theta,h,c)$ from \eqref{SheppResult0}. To prove \eqref{BMDensXThHCR}, we  multiply \eqref{BMDensR} by $p(\theta,h,c)$.
\qed

Using \eqref{CHAleft} and \eqref{mnMeander2}, we calculate the moments of the conditional Brownian motion:
\begin{theorem} \label{BrownMom}
Under the assumptions of Theorem \ref{CHAThm}, the expectation of a Brownian motion satisfies
$E[B(t) | c, h, \theta]$. 
%$Var[B(t) | c, h, \theta]$:
\BEQ \label{MnCThHa}
E[B(t) | c, \theta,h]= h - \sqrt{\theta} M_1 (1 -\frac{t}{\theta},\frac{h}{\sqrt{\theta}}) \ = %\ {\rm for}\ \  t \le \theta
\NEQ 
\BEQ \label{MnCThHaa}
h - \sqrt{\theta} [\frac{t}{h\sqrt{\theta}} + \frac{h(\theta-t)}{\theta^2}] {\rm erf}(\frac{h\sqrt{\theta - t}}{\sqrt{2t\theta}}) +\frac{\sqrt{2t(\theta-t)}}{\sqrt{\pi \theta}} \exp(\frac{-(\theta-t)h^2}{2t\theta}) 
\NEQ 
for $t \le \theta$. 
For $t \ge \theta$, the expectation satisfies:
\BEQ\label{MnCThHb}
E[B(t) | c, \theta,h]= h - \sqrt{1-\theta} M_1(\frac{t-\theta}{1-\theta},\frac{h-c}{\sqrt{1-\theta}}) \ \ {\rm for}\ \ t \ge \theta \ \ .   % {\rm where}\ \ q \equiv (H-c) / \sqrt{1-\theta}.  
\NEQ
For $t\le \theta$,  $E[B^2(t) | c, \theta,h]= h^2 - 2h\sqrt{\theta}M_1 (s,\frac{h}{\sqrt{\theta}}) + \theta M_2(s,\frac{h}{\sqrt{\theta}})$. Thus the variance satisfies
\BEQ \label{VrCThHa}
Var[B(t) | c, \theta,h]= \theta Vr (1 -\frac{t}{\theta},\frac{h}{\sqrt{\theta}}) \ \ {\rm for}\ \  t \le \theta \ ,
%h - 3t(\theta-t) + \sqrt{\theta} M_1 (1 -t/\theta) \ \ {\rm for}\ \  t \le \theta
\NEQ 
\BEQ \label{VrCThHb}
Var[B(t) | c, \theta,h]= (1-\theta) Vr(\frac{t-\theta}{1-\theta},\frac{h-c}{\sqrt{1-\theta}}) \ \ {\rm for}\ \ t \ge \theta \ \ ,
\NEQ
where $Vr(s,c)\equiv  M_2(s,c)-M_1(s,c)^2=3s(1-s)+ c^2s^2 - M_1(s,c)^2$ .
\end{theorem}

\Prf In the representations, \eqref{CHAleft} and \eqref{CHAright}, the conditional density of the Brownian motion
(on either side of the $\argmax$, $\theta$), is a linear transformation of the Brownian meander density. Thus the moments of the conditional density, \eqref{BMDensR} and \eqref{BMDensL}, are linear transformations
of the Brownian meander moments given in Theorem \ref{BMmom}. The representations in \eqref{MnCThHa} - \eqref{VrCThHb}
follow by substituting the appropriate transformations of the Brownian meander moments into  \eqref{CHAleft} - \eqref{CHAright}.

\section{Moments Given High and Its Location} \label{EHThsect}

To calculate the expectation of $B(t)$ given only its maximum and the first location of its maximum, we
integrate the expressions in \eqref{BMDensL} and \eqref{BMDensR} over the correct conditional probability density
given by \eqref{SheppResult0} and \eqref{eq:denThH}:
\BEQ \label{EHThH}
p(x,t|\theta,h)= \int_{-\infty}^h p(x,t | c, \theta,h) p(\theta, h, c) dc /p(\theta,h) \ .
\NEQ

\begin{theorem} \label{denTH}
The conditional density, $p(B(s)=x|\theta,h)$ satisfies
\begin{align} \label{BMDensXThH}
  p(B(s)=x |\theta,h)= \frac{(h-x)}{\pi(1- \theta)^{\frac{1}{2}} (\theta-s)^{1.5}}  g_{s}(h-x,h)
\exp\left(\frac{-(h-x)^2}{2(\theta-s)}  \right) 
\end{align}
for $t<\theta$.
Here $g_{1-s}(x,y) = \phi_{1-s}(x-y) - \phi_{1-s}(x+y)$.
For $t>\theta$, the conditional density is
\begin{align} \label{BMDensXThHR}
  p(B(s)=x|\theta,h)= \frac{h(h-x)}{\pi \theta^{1.5} (s-\theta)^{1.5}}  \erf(\frac{h-x}{\sqrt{2(1-s)}})
  \exp(\frac{-(h-x)^2}{2(s-\theta)}-\frac{h^2}{2\theta}) \ \ .
\end{align}
% integrate (\exp(-.57 (c-x)^2/2)- \exp(-.57(c+x)^2/2))dx from 0 to infty
\end{theorem}

\Prf 
We integrate \eqref{BMDensXThHC} and \eqref{BMDensXThHCR} with respect to the $close$, $c$ and then divide by $p(\theta, h)$ as given by \eqref{eq:denThH}. We use $\int_{-\infty}^h  g_{1-s}(z,h-c) dc$ $= \int_{0}^{\infty} \phi_{1-s}(x-z) - \phi_{1-s}(x+z) dx = \erf(\frac{z}{\sqrt{2(1-s)}})$.
\qed

To calculate the expectation of $B(t)$ given only its maximum and the first location of its maximum, we need only
integrate the expressions in \eqref{mnMeander1} - \eqref{MnCThHb} over the same conditional probability density:
\BEQ \label{EHTh1}
E[B(t|\theta,h)]= \int_{-\infty}^h E[B(t) | c, \theta,h)] p(\theta, h, c) dc /p(\theta,h) 
\NEQ
where $p(\theta, h, c)$ is given by \eqref{SheppResult0} and $p(\theta,h)$ is given by \eqref{eq:denThH}.
Thus $p(\theta, h, c)/p(\theta,h) =\frac{(h-c)}{(1-\theta)} \exp(\frac{-(h-c)^2}{2(1-\theta)})$.  We change variables
from $c$ to $x\equiv \frac{(h-c)}{\sqrt{1-\theta}}$ and apply \eqref{AM1Int} - \eqref{AM1IntC} the Appendix A.
\begin{theorem}
For $t>\theta$, the expectation becomes
\BEQ \label{EHTh2}
E[B(t|\theta,h)]= h- \sqrt{1-\theta}\int_0^{\infty} x M_1(s,x) e^{-x^2/2} dx = h -\sqrt{1-\theta}G_{11}(s) \ \ ,
\NEQ
where $s \equiv  \frac{t-\theta}{1-\theta}$ and $G_{11}(s)$ is defined as
\BEQ \label{G11Def} %{EHTh5}
G_{11}(s) = \sqrt{\frac{2}{\pi}}\left[\tan^{-1}(\sqrt{\frac{s}{1-s}}) + \sqrt{s(1-s)}\right] \ \ .
%\ \ {\rm for}\ \ t \ge \theta \
\NEQ
%where $s \equiv  \frac{t-\theta}{1-\theta}$ and $G_{11}(s)$ is defined by \eqref{EHTh5}. % and $\kappa = s/(1-s)$. %We use the identities:
%$\int_0^{\infty} \erf(\frac{ax}{\sqrt{2}}) \exp(\frac{-x^2}{2})dx =\sqrt{\frac{2}{\pi}} \tan^{-1}(a)$ and
%$\sqrt{\frac{\pi}{2}} \int_0^{\infty}x^2 \erf(\frac{ax}{\sqrt{2}}) \exp(\frac{-x^2}{2})dx =  \frac{a}{a^2 +1} + \tan^{-1}(a)$.
%Note that $E[B(t|\theta,h)] = h +g(t,\theta)$ where $g$ is defined by \eqref{EHTh5}.
Thus $E[B(t|\theta,h)] -h $  is independent of $h$ for $t \ge \theta$.
For $t \le \theta$, $E[B(t) | c, \theta,h]$ is independent of $c$ and satisfies
$E[B(t) |\theta,h]=E[B(t) | c, \theta,h]$ as given by \eqref{MnCThHa} - \eqref{MnCThHaa}.
For $t \le \theta$, the variance is independent of the close and is given by \eqref{VrCThHa}:
\BEQ \label{VHTh1L}
Var[B(t|\theta,h)]= \theta Vr (1 -\frac{t}{\theta},\frac{h}{\sqrt{\theta}}) \ \ {\rm for}\ \  t \le \theta \ \ .
\NEQ
For $t \ge \theta$, the variance of $B(t|\theta,h)$
\BEQ \label{VHTh}
Var[B(t|\theta,h)]=(1-\theta) [3 s -s^2]  - (1-\theta)G_{11}(s)^2  \ \ .
%\ {\rm for}\ \ \ t \ge \theta 
\NEQ
\end{theorem}

\Prf The proof of \eqref{EHTh2} is in \eqref{AM1Int} - \eqref{AM1IntC}. The variance of $B(t|\theta,h)$ is calculated by the ensemble average of $E[B^2(t|c,\theta,h)]$ and subtracting $E[B(t|\theta,h)]^2$. For $t \le \theta$, $E[B(t)^2 | c, \theta,h]$ is also independent of $c$. This proves \eqref{VHTh1L}.
For $t \ge \theta$,
\BEQ \label{VHThP}
Var[B(t|\theta,h)]= \int_{-\infty}^h E[B^2(t) |c,\theta,h]p(c|\theta, h) dc -  E[B(t|\theta,h)]^2
\NEQ
\BEQ \label{VHTh2R}
= h^2-2\sqrt{(1-\theta)}hE[B(t|\theta,h)] + (1-\theta)\int_0^{\infty} x M_2(s,x) e^{\frac{-x^2}{2}} dx - E[B(t|\theta,h)]^2 \ \ .
\NEQ
To evaluate \eqref{VHTh2R}, we insert \eqref{AM2Int} and \eqref{EHTh2} into \eqref{VHTh2R}.
\qed

%For $t \ge \theta$, we evaluate the conditional variance of $B(t)$ using \eqref{AM2Int}:
By ignoring the $c$ dependence in $B(t|c,\theta,h)$, the variance increases by
$\int_{-\infty}^h E[B(t|c,\theta,h)]^2$ $p(c|\theta, h) dc -$ $E[B(t,\theta,h)]^2$. One can think of this as an Analysis of Variance term.
The increase occurs because the final value, $c$, is now a random variable.
% stochastic, this analysis of variance has the form: 

Similar but more complicated moment calculations are possible for the case of the previous section where $(\theta,c)$ are given and one integrates over $h$. In \cite{Riedel}, analytic expressions for $E[B(t|c,h)]$ and $Var[B(t|c,h)]$ are given. It may be
possible to rederive these results by integrating $p(x,t | c, \theta,h) p(\theta, h, c)$ $d\theta /p(c,h)$. This integration appears to be untractable.

\section{Mean and Variance Given Only Argmax \label{EVAsect}}

Now consider $B(t|{\rm \argmax}\ B =\theta)$.
To compute the mean and variance of $B(t|\theta)$, we integrate \eqref{EHTh1} - \eqref{EHTh2} with respect to $p(h|\theta)dh$
where $p(h|\theta)$ is given by \eqref{eq:condenThH}:
\begin{theorem} \label{argmaxThm}
 For $t<\theta$,
\BEQ \label{ETh1}
E[B(t|\theta)]= \int_{0}^{\infty} E[B(t) |\theta,h] p(h|\theta) dh = \sqrt{\frac{\pi\theta}{2}} - \int_{0}^{\infty} M_1 (1 -\frac{t}{\theta},\frac{h}{\sqrt{\theta}}) \frac{h}{\sqrt{\theta}} e^{\frac{-h^2}{2\theta}}dh
\NEQ
\BEQ \label{ETh1b}
=\sqrt{\frac{\pi \theta}{2}}-\sqrt{\theta}G_{11}(s) 
= \sqrt{\frac{\pi \theta}{2}}-\sqrt{\frac{2\theta}{\pi}}\left[\tan^{-1}(\sqrt{\frac{s}{1-s}}) + \sqrt{s(1-s)}\right] \ \ ,
\NEQ
where $s\equiv 1 -\frac{t}{\theta}$ and $G_{11}$ is given in \eqref{G11Def}.
For $t>\theta$, the expectation becomes
\BEQ \label{ETh5}
E[B(t|\theta)]= \sqrt{\frac{\pi \theta}{2}} - \sqrt{\frac{2(1-\theta)}{\pi}}\left[\tan^{-1}(\sqrt{\frac{s}{1-s}}) + \sqrt{s(1-s)}\right] \ \ {\rm for}\ \ t \ge \theta \
\NEQ
where $s \equiv  \frac{t-\theta}{1-\theta}$. % and $\kappa = s/(1-s)$.

\noindent For the variance, we first evaluate $E[B^2(t|\theta)]$ and then subtract off $E[B(t|\theta)]^2$.
%$h^2 - 2\sqrt{\theta}hM_1 (s,\frac{h}{\sqrt{\theta}}) + \theta Me_2(s,\frac{h}{\sqrt{\theta}})$
\BEQ \label{B2Th1La}
E[B^2(t|\theta)]= 2 \theta- 4 \theta \sqrt{s} +\theta [3 s -s^2] \ \ \  {\rm for}\ \ \  t \le \theta
\NEQ
%Then $Var[B(t|\theta)]=E[B^2(t|\theta)]-E[B(t|\theta)]^2$.
For $t \ge \theta$, the variance simplifies:
\BEQ \label{VrThR}
Var[B(t|\theta)]= Var[h|\theta] + (1-\theta) \left[3s -s^2 - G_{11}(s)^2 \right] \ \ ,
\NEQ
where $Var[h|\theta]= (2-\pi/2) \theta$.
\end{theorem}

\Prf We use the integrals in  \eqref{AM1Int} - \eqref{AM1IntC} of the Appendix A and set $s= 1- t/\theta$. 
To prove \eqref{B2Th1La}, note
\BEQ \label{B2Th1L}
E[B^2(t|\theta)]= \int_0^{\infty} E[B^2(t) |\theta,h] p(h|\theta)dh
=\int_0^{\infty}\left[\theta x^2 - 2\theta x M_1 (s,x) + M_2(s,x)  \right] x e^{\frac{-x^2}{2}}dx
\NEQ
and insert the expression in \eqref{AM1TInt} - \eqref{AM1TIntB}.

%$M_2(s,t,c)=\frac{3s(h-s\sqrt{\theta})}{h}+ \theta s^2/h^2$.
%\BEQ \label{B2Th2L}
%E[B^2(t|\theta)]= \theta\int_0^{\infty} [\frac{3s(h-s\sqrt{\theta})}{h}+ \frac{\theta s^2}{h^2} ]h\exp(-h^2/2\theta) dh = 
%\theta [3 s -s^2]
%\NEQ
%where $s\equiv 1 -\frac{t}{\theta}$. Similarly,
%\BEQ \label{VTh1R}
%E[B^2(t|\theta)] =  (1-\theta) [3 s -s^2] \ \ {\rm for}\ \ t \ge \theta \ .
%\NEQ
%We note that $E[B(t|\theta,h)]$ 
Our expression matches the zero dimensional results derived in the second section, $E[h|\theta] = \sqrt{\pi \theta/ 2}$, $Var[h|\theta] = (2-\pi/2) \theta$.
%Figure %\ref{fi:AHMeanCombo} and Figure \ref{fi:AHVarCombo}
%compare our simulations with the analytic expressions in \eqref{ETh1} - \eqref{ETh5}

\begin{figure}[htbp]
  \centering
\includegraphics[scale=\FSCL]{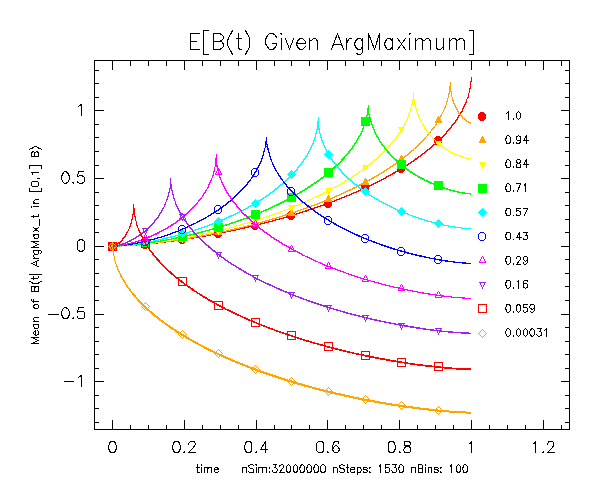}  
\caption{Expectation of $B(t)$ given $\theta\equiv \argmax\{B(t)\}$ for various values of $\theta$.
Each color has two curves, a theoretical curve from Theorem \ref{argmaxThm} and the mean value of the simulation for the given parameter bin. \n
  To compute the simulation expectation,
  we use an ensemble of 32,000,000 realizations computed with 1500 steps and bin the results into 100 bins in $\theta$ space.
  The values of $\theta$ for each curve are given in the legend.}
  \label{fi:MeanArrGivenArgMax}
\end{figure}

 Let $g(t,\theta)=E[B(t|{\rm \argmax}\ B =\theta)]$. %and $var(t,\theta)\equiv Var[B(t)|{\rm \argmax}\ B= \theta]$.
Figure \ref{fi:MeanArrGivenArgMax} plots $g(t,\theta)$ for ten quantile values of the bins in $\theta$ ranging from the second smallest bin to the second largest bin.
Unsurprisingly, the $g(t,\theta)$ is maximized at $t=\theta$ for
fixed $\theta$. The curves look monotone for $\theta$ very near $0$ or $1$, but this is a graph and simulation resolution effect.
Even at $nbins=160$, there is still a small amount of artificial broadening of the peak due to averaging  $g(t,\theta)$ over the bin. Using the density \eqref{eq:denHC0}, we calculate $E[h|\theta] = \sqrt{\pi \theta/2}$. 

%A naive calculation of $\argmax$ on a discrete set of $nSteps$ timepoints results a very small number of discrete values of $\theta$. Having many fewer $\argmax$ values than our large number of simulations, $nSim$, can play havoc with the adaptive quantile gridding. To prevent this, we use a three point parabolic interpolation around the actual $\argmax$ of the finite time grid. If the actual $\argmax$ occurs at the left most or rightmost timepoint, we take a random point in the time delta about the first or last point. It is possible to more accurately estimate the expectation of arg-maximum of a Brownian process given two or three data points near the maximum \cite{AC,Devroye,McLeish}  

\begin{figure}[htbp]
\centering
\includegraphics[scale=\FSCL]{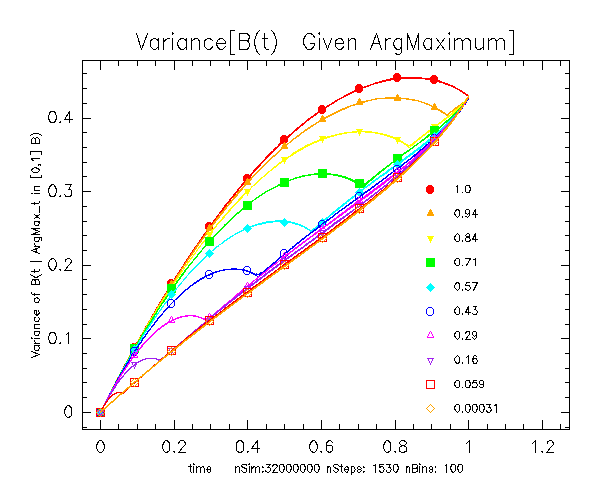}
\caption{Variance of $B(t)$ given $\theta\equiv \argmax\{B(t)\}$ for various values of $\theta$. When $\theta$ is small, the variance grows
roughly linearly in time. For larger values of $\theta$, the variance grows faster initially and then converges to the linear value.}
\label{fi:VarArrGivenArgMax}
\end{figure}

%Figure \ref{fi:VarArrGivenArgMax} plots $Var[B(t)|{\rm \argmax}\ B= \theta]$ for ten quantile values of the bins in $\theta$.
Figure \ref{fi:VarArrGivenArgMax} plots the variance, $var(t,\theta)=Var[B(t)|{\rm \argmax}\ B= \theta]$.
For small values of $\theta$,  $var(t,\theta)$ grows linearly for most time,
then it accelerates as $\theta $ approaches $1$. For larger values of $\theta$, it appears that $\frac{\partial var(t,\theta) }{\partial t}$
is nearly a constant for $t$ near zero. The growth rate decelerates in time and this deceleration occurs earlier for smaller $\theta$.
%We note that the value of $var(t=1,\theta) = Var[\sqrt(2E)] = 2*(1 -\pi/4) \approx 0.4292$.
From \eqref{eq:B1condTh}, $Var[ B(1) | \theta] = 2(1 - \pi/4) \approx 0.4292$.

\section{Overview of Simulation Results \label{SimDesc}}
%\graphicspath{{/home/kriedel/brownOut/pubAH/}}
%\graphicspath{{/home/kriedel/brownOut/pubAH/},{/home/kriedel/brownOut/pubCAHo/}}

We generate a large number of Brownian paths, bin the paths in $(close, \max, \argmax)$ space
and calculate the mean and variance for each time and bin. We order the coordinates of phase
space, $(q_1,q_2, q_3)$, so that $q_1 = B(1)$, $q_2=\ \max_{0 \le t \le 1}B(t)$ and $q_1=\ \argmax_{0 \le t \le 1}B(t)$.
To see  $E[B(t) |givens]$ and $Var[B(t) |givens]$, we simulate 10 to 36 million Brownian motions. 
A detailed description of the simulation method may be found in Appendix B and \cite{Riedel}.
We display simulations using $1500$-$1850$ time steps for each realization.

%\graphicspath{{/home/kriedel/brownOut/pubCAH/},{/home/kriedel/brownOut/pubCHA/},{/home/kriedel/brownOut/pub2/}}

For each bin, we compute the expectation and variance using \eqref{MnCThHa} - \eqref{VrCThHb}.
In Appendices C, D and E, we present a number of figures displaying and explaining the behavior of
$E[B(t) | c, \theta,h]$ and $Var[B(t) | c, \theta,h]$ as $t,c,\theta,h$ vary. %We omit these for space reasons.

We evaluate the discrepancy/mean square error (MSE) between the theoretical formulae and our simulations. For each bin in parameter space, we compare the difference between \eqref{MnCThHa} and the average
of the curves in the bin. We then compute the mean squared error averaged over time:
\BEQ
MSE(c,\theta,h)=\int_0^1\left(E[B(t|c, \theta, h)] -  Avg_j[\tilde{B}_j(t|c, \theta, h)]\right)^2 dt
\NEQ
where $\tilde{B}_j$ is the $j$ curve in the bin and $Avg_j$ is the average over all of the paths in the given bin in parameter space.
We then sort the bins from worst to best MSE. Not surprisingly
the worst fitting bins often occur where the bins are the largest and the bias error, from having curves
with slightly different parameters, is the largest. 

\begin{figure}[htbp]
\centering
\includegraphics[scale=\FSCL]{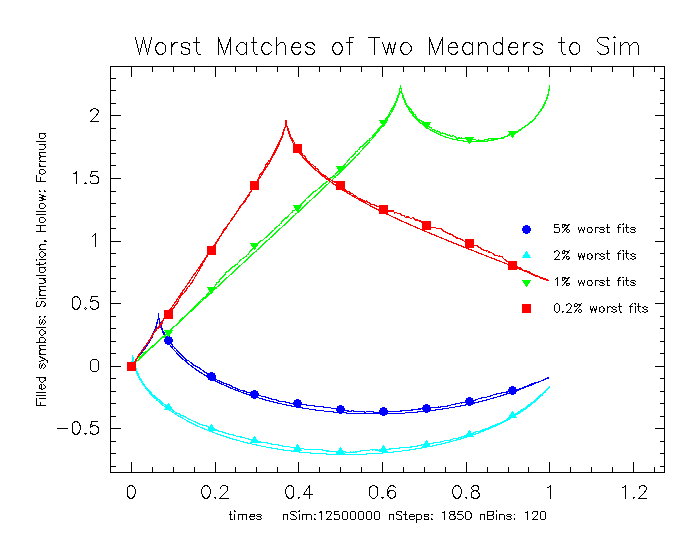}
\caption{Comparison of  Formula and Simulation: Each color has two curves, a theoretical curve
  from Theorem \ref{BrownMom} and the mean value of the simulation for the given parameter bin.
  The curves overstruck by symbols are the simulation curves. The analytic formula curves have the same color but no symbol.\n
Blue:  $5\%$  worst MSE Mean:0.00039 at close:-0.0891, high:0.431, argmax:0.065\n
Cyan: $2\%$ worst MSE Mean:0.000507 at close:-0.166, high:0.101, argmax:0.00374\n
Green: $1\%$ worst MSE Mean:0.000608 at close:2.246, high:2.27, argmax:0.643\n
Red: $0.2\%$ worst MSE Mean:0.00091 at close:0.682, high:1.991, argmax:0.37
 }
\label{fi:MeanderMeanCombo}
\end{figure}

Overall, the concatenated meander expectation fits the simulation expectation very well. To further illustrate this, we overlay all four sets of curves on the same plot. Figure \ref{fi:MeanderMeanCombo} plots the curves that correspond to
the worst $5\%$, $2\%$, $1\%$ and $0.2\%$ difference between the simulation conditional mean and \eqref{MnCThHa} - \eqref{MnCThHb}.
There are $120^3 \ge 1.7$ million bins to compare.
The curves overstruck by symbols are the simulation curves. The analytic formula curves have the same color but no symbol.

\begin{figure}[htbp]
\centering
\includegraphics[scale=\FSCL]{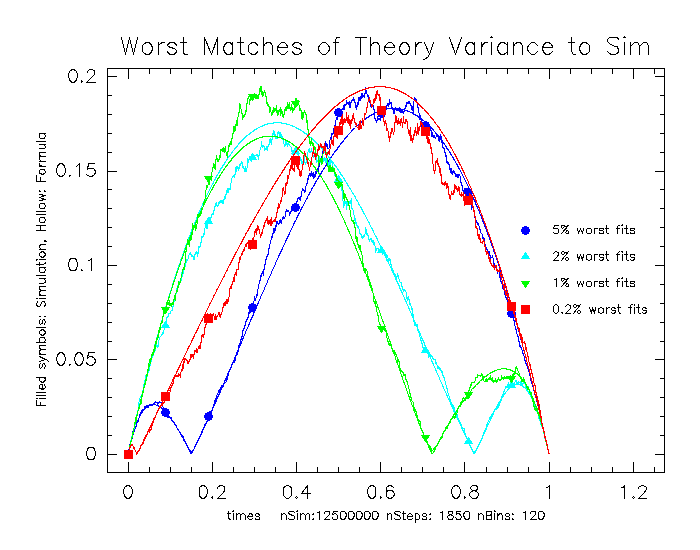}
\caption{Comparison of Meander Variance Formula and Simulation: Each color has two curves, \n
  a theoretical variance  and the variance of the simulation for the given parameter bin. \n
The symbols are on the curve from the simulation. \n
Blue: $5\%$ worst MSE Var:.000054 at close:-1.743, high:0.55, argmax:0.15 \n
Cyan:$2\%$ worst MSE Var:.0000775 at close:1.246, high:2.177, argmax:0.822 \n
Green:$1\%$ worst MSE Var:.0000981 at close:2.531, high:3.164, argmax:0.723 \n  
Red:$0.2\%$ worst MSE Var:0.000159 at close:-1.437, high:0.461, argmax:0.0208  }
\label{fi:MeanderVarCombo}
\end{figure}

We conclude with a comparison of the variance of the simulation with the concatenated meander in \eqref{VrCThHa} - \eqref{VrCThHb}.
Figure \ref{fi:MeanderVarCombo} plots the curves that correspond to
the worst $5\%$, $2\%$, $1\%$ and $0.2\%$ difference between the simulation conditional variance and \eqref{VrCThHa} - \eqref{VrCThHb}.
In Appendix C, we display the analogous plots for $E[B(t|\theta, h)$. As expected,
averaging over values of $c$ reduces the noise in the simulations and reduces the discrepancy.
  
%\vspace{-5mm}

\section{Summary}

\begin{figure}[htbp]
  \centering
\includegraphics[scale=\FSCL]{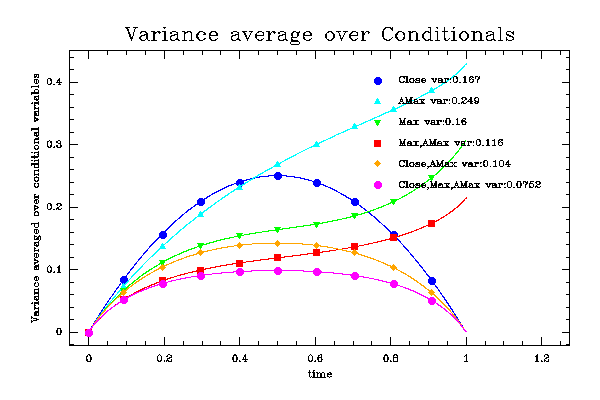}  
\caption{Comparison of $Var[B(t|close)]$,$Var[B(t|close)]$, versus $Var[B(t|close, \argmax)]$ versus $Var[B(t|close, \argmax, high)]$.
  These curves are the ensemble average over all of the simulations.}
\label{fi:CAH:AvgVarGivenCAH}
\end{figure}

We have investigated the distribution, expectation and variance of Brownian motion conditional on the values of its maximum and its location as well
as its final value. We give formulae for the mean and variance of a Brownian meander. Using the Williams construction, we evaluate the expectation and variance, $E[B(t) | c, \theta,h]$ and $Var[B(t) | c, \theta,h]$.
To analyze $B(t|\theta,h)$ and $B(t|\theta)$, we integrate the distribution and moments with respect to the close, $c$, and then with
respect to the maximum, $h$. Similar but more complex formulae can be calculated for $B(t|\theta,c)$.
Our simulations show good agreements with the analytic expressions. Many interesting features are displayed in the figures of the mean and variance.
Figure \ref{fi:CAH:AvgVarGivenCAH} displays the ensemble average of the variance of $B(t|givens)$. For example,
\BEQ \label{Vens}
V(t|{\mathrm givens}=(c,\theta,h))\equiv    \int Var[B(t|c,\theta,h)] dp(c,\theta,h) \ .\
%V(t|h,c)\equiv \int Var[B(t|h,c)] dp(h,c) \ .
\NEQ
We compare the ensemble averages of
$Var[B(t|\argmax)]$ with $Var[B(t|close)]$,$Var[B(t|close, \argmax)]$ and $Var[B(t|close, \argmax, high)]$. In all case, we compute the variance for each value of the givens separately and then average over the distribution of the givens.
From Figure \ref{fi:CAH:AvgVarGivenCAH}, we see that $\argmax$ is a much less valuable statistic than either the $close$ or $max$.
The statistics, $(close, \argmax)$, are more valuable than the pair of statistics $(\max,\argmax)$. If the final value ($close$)
is known, the variance is symmetric in time. Otherwise, the variance increases in time.
%In the next section, we  display $E[B(t |close, \argmax)]$  and $Var[B(t |close, \argmax)]$.
%Section \ref{EVAsect} does the same for $B(t|\argmax)$ and Section \ref{EVCAH} does the same for $B(t|close, \argmax, max)$. 
%%In the section \ref{ECAH}, we add in the maximum value of $B(t)$.
%%begin our gallery of the moments of Brownian motion by

We conclude by merging our results with those  in \cite{Riedel}. To measure how accurately we can estimate $B(t|\theta,h,c)$, we consider the time integral
of the conditional variance averaged over the given variables:  $\int \int_0^1Var[B(t|c,\theta,h)]dt dp(c,\theta,h)$. If only the final value, $c$,
is specified, the value of the integral is $1/6$. Table 1 presents the corresponding values for our simulations in this paper.
\begin{table}[ht] 
%\caption{Time Averaged Variance by Givens}
% title of Table
\centering 
%\hline
\begin{tabular}{|c|c|c|}
\hline 
{\bf Givens}&{\bf Var}& {\bf Var*6} \\
\hline
Start point only& 1/2 & 3\\
Close & 1/6 &1\\
High &0.1602  & .9612\\
ArgMax & .2487 & 1.492\\
Close, High & 0.0990&.5938\\ %.0989651305985
Close, ArgMax & 0.1037 &.6222\\
Argmax, High & 0.11585 &0.6951 \\
%High, Low & 0.09911& .5947\\
%ArgMax, ArgMin & 0.1574 & 0.9444\\
Close, High, ArgMax &  0.07535&.4521\\ 
Close, High, Low & 0.0701 &.4206\\  %0.0701043061526
\hline \hline
\end{tabular}
\label{EnsVarTab} %\label{tabVar}
%\caption
\caption{Expected time average variance reduction. We multiply the variance by $6$ in the third column to compare with knowing only the final value, $c$.}
\end{table}
This shows that the statistics triple $(Close, High, Low)$ is more valuable than the triple $(Close, High, ArgMax)$.
In all cases, adding addition information such as the value of the high and/or its location significantly
reduces the time averaged variance in comparison with that of the Brownian bridge (i.e only using the final value of $B(t=1)$).

In our sister article \cite{Riedel}, we show using the high, low and close substantially improves the
estimation of the log return of the SP500 over using just the close. We can reproduce this
calculation/application in this article, but it will be artificial. The reason is that
the location of the maximum is seldom if ever used to parameterize a stochastic process.
In contrast, chartist analysis in finance heavily uses the values of the high and low.
Note that Table 1 shows that using the value of the maximum is better than using the value of arg maximum.

\section{Appendix A: Meander Moments}
  We now compute moments
  \begin{align} \label{meanderMom}
M_k(s,t,c) &= \int_0^{\infty} x^k P(B^{me}(s)=x|B^{me}(t)=c) dx \\    
  &= \int_0^{\infty}\left( \phi_{t-s}(c-x) -\phi_{t-s}(c+x)\right) \frac{x^k t^{3/2}}{cs^{3/2}} \exp(-x^2/2s) \exp(c^2/2t)  dx \\
&= \frac{t^{3/2}}{cs^{3/2}}\frac{2}{\sqrt{2\pi \tau}} \exp(-c^2/2h) \exp(c^2/2t)  \int_0^{\infty} x^k \exp(-a x^2)\sinh(bx) dx \ \ ,
\end{align}  
where $\tau \equiv t-s$, $a \equiv\frac{1}{2\tau}+ \frac{1}{2s} = \frac{t}{2s(t-s)}$ and $b\equiv \frac{c}{t-s}$.
Note $M_0(s,t,c)=1$: 
\BEQ  
M_0(s,t,c) = \frac{t^{3/2}}{cs^{3/2}}\frac{1}{\sqrt{2\pi \tau}} \exp(-c^2/2h) \exp(c^2/2t)  \frac{\sqrt{\pi} b }{2\sqrt{a^3}} e^{\frac{b^2}{4 a}}  =1 \ .
\NEQ  
Here we use $b^2/4a =c^2 s/2t(t-s)$ and $\frac{1}{2h}- \frac{1}{2t} = \frac{s}{2t(t-s)}$ to show that the exponential terms cancel. %\ref(mnMeander4}.
%Similarly, $M_1(s,t,c)$ satisfies
The first moment of the meander satisfies
\begin{align} 
M_1(s,t,c) &= \frac{t^{3/2}}{cs^{3/2}}\frac{1}{\sqrt{2\pi \tau}} \exp(-c^2/2h) \exp(c^2/2t)\left[ \frac{\sqrt{\pi}[2a+ b^2] e^{\frac{b^2}{4 a}}{\rm erf}(\frac{b}{2 \sqrt{a}})}{4a^{2.5}} +\frac{b}{2a^2}\right] \\
  &= \frac{[t-s+ sc^2/t] {\rm erf}(\frac{c\sqrt{s}}{\sqrt{2t(t-s)}})}{c} +\frac{\sqrt{2s(t-s)}}{\sqrt{\pi t}} \exp(-sc^2/2t(t-s)) \ \ .
\end{align}
%where $a=\frac{1}{2h}+ \frac{1}{2s} = \frac{t}{2s(t-s)}$ and $b= \frac{c}{t-s}$. Note that $b^2/4a =c^2 s/2t(t-s)$,
%where $b/2a^2 =2 c s^2(t-s)/t^2$.
Similarly, the second moment of the meander satisfies
\begin{align} \label{mnMeander2C}
M_2(s,t,c) &= \frac{t^{3/2}}{cs^{3/2}}\frac{1}{\sqrt{2\pi \tau}} \exp(-c^2/2h) \exp(c^2/2t)  \frac{\sqrt{\pi} b (6 a + b^2)e^{\frac{b^2}{4 a}} }{8a^{7/2}} = \\
  &= \frac{t^{3/2}}{cs^{3/2}}\frac{1}{\sqrt{2(t-s)}}  \frac{ b (6 a + b^2) }{8a^{7/2}} 
  = \frac{t^{3/2}}{s^{3/2}}\frac{1}{\sqrt{2(t-s)^{1.5}}}  \frac{ (1.5 + b^2/4a) }{2a^{5/2}} \\
  %&= \frac{t^{3/2}}{s^{3/2}}\frac{1}{\sqrt{8(t-s)^{1.5}}}  (1.5 + c^2 s/2t(t-s)) \left( \frac{2s(t-s)}{t} \right)^{5/2} \\
&= (3 + c^2 s/t(t-s)) \left( \frac{s(t-s)}{t} \right)     = \frac{3s(t-s)}{t}+ c^2s^2/t^2 \ .
\end{align}  
%where $a=\frac{1}{2h}+ \frac{1}{2s} = \frac{t}{2s(t-s)}$ and $b= \frac{c}{t-s}$. Note that $b^2/4a =c^2 s/2t(t-s)$.
The integrations follow from
%$\int_0^{\infty} x \exp(-a x^2 \pm bx) dx=  \frac{1}{2a} +
%  \frac{\sqrt{\pi} b }{4\sqrt{a^3}} e^{\frac{b^2}{4 a}} \left[{\rm erf}(\frac{b}{2 \sqrt{a}})\pm 1 \right]$
\begin{align}  
\int_0^{\infty} x \exp(-a x^2)\sinh(bx) dx &=   \frac{\sqrt{\pi} b }{4\sqrt{a^3}} e^{\frac{b^2}{4 a}} \\ 
 \int_0^{\infty} x^2 \exp(-a x^2)\sinh(bx) dx &=
\frac{\sqrt{\pi}[2a+ b^2] e^{\frac{b^2}{4 a}}{\rm erf}(\frac{b}{2\sqrt{a}})}{8a^{2.5}} +\frac{b}{4a^2} \\
\int_0^{\infty} x^3 \exp(-a x^2)\sinh(bx) dx &= \frac{\sqrt{\pi} b (6 a + b^2)e^{\frac{b^2}{4 a}} }{16a^{7/2}} \quad \ \ .
\end{align}

%$\int_0^{\infty} x^3 \exp(-a x^2 \pm bx) dx= \frac{\sqrt{\pi} b (6 a + b^2)e^{\frac{b^2}{4 a}} \left[{\rm erf}(\frac{b}{2 \sqrt{a}}) \pm 1 \right]}{16 a^{7/2}} + \frac{4 a + b^2}{8a^3}$

%$\int_0^{\infty} x^2 \exp(-a x^2 \pm bx) dx=
%  \frac{\sqrt{\pi}[2a+ b^2] e^{\frac{b^2}{4 a}}\left[1 \pm {\rm erf}(\frac{b}{2 \sqrt{a}}) \right]}{8a^{2.5}} \pm \frac{b}{4a^2}$

%$\int_0^{\infty} \exp(-x^2/2)\exp(-b (x-a)^2/2) dx= \sqrt{\frac{\pi}{2}}  \frac{e^{-\frac{a^2b}{2b+2}}\left[1+ \erf(\frac{ab}{\sqrt{2b+2}})\right]}{\sqrt{b+1}}  $.
For Section \ref{EHThsect}, we need the following integral:
\begin{align}
\int_0^{\infty} x M_1(s,x) e^{\frac{-x^2}{2}} dx &=
\int_0^{\infty} (1-s+sx^2) e^{\frac{-x^2}{2}}\erf(x\sqrt{\frac{\kappa}{2}})  +\sqrt{\frac{2}{\pi}}\sigma(s) xe^{\frac{-x^2}{2(1-s)}} dx \label{AM1Int} \\
%\NEQ \BEQ
&=\sqrt{\frac{2}{\pi}}\left[\tan^{-1}(\sqrt{\kappa}) + s \frac{\sqrt{\kappa}}{\kappa +1}  + \sqrt{s(1-s)^3}\right]  \label{AM1IntB}\\
%e^{\frac{-h^2}{2(1-\theta)}}
%h- \sqrt{\frac{2(1-\theta)}{\pi}}\left[\tan^{-1}(\sqrt{\kappa}) + s \frac{\sqrt{\kappa}}{\kappa +1} \right] + \sqrt{\frac{2s(1-s)^3(1-\theta)}{\pi}} %e^{\frac{-h^2}{2(1-\theta)}}
&= G_{11}(s)\equiv\sqrt{\frac{2}{\pi}}\left[\tan^{-1}(\sqrt{\frac{s}{1-s}}) + \sqrt{s(1-s)}\right] \ ,
\label{AM1IntC}
\end{align}
where $\kappa = s/(1-s)$ and $\sigma(s)\equiv \sqrt{s(1-s)}$. We use the identities:
$\int_0^{\infty} \erf(\frac{ax}{\sqrt{2}}) exp(\frac{-x^2}{2})dx =\sqrt{\frac{2}{\pi}} \tan^{-1}(a)$ and
%$\int_0^{\infty}x \erf(\frac{ax}{\sqrt{2}}) \exp(\frac{-x^2}{2})dx =\frac{a}{\sqrt{a^2 +1}}$
$\sqrt{\frac{\pi}{2}} \int_0^{\infty}x^2 \erf(\frac{ax}{\sqrt{2}}) \exp(\frac{-x^2}{2})dx =  \frac{a}{a^2 +1} + \tan^{-1}(a)$.
To evaluate the variance, we need
\BEQ \label{AM2Int}
\int_0^{\infty} M_2(s,x,1)x e^{\frac{-x^2}{2}}dx  = \int_0^{\infty}  (3s(1-s)+ x^2s^2) x e^{\frac{-x^2}{2}}dx=
[3 s(1-s)+2 s^2] = [3 s -s^2] \ .
\NEQ

To evaluate  $E[B^2(t|\theta)]$ in \eqref{B2Th1La}, note
%For Section \ref{EVAsect}, we define $G_{1,2}(s)$ as 
\BEQ \label{AM1TInt}
\int_0^{\infty} x^2 M_1(s,x) e^{\frac{-x^2}{2}} dx =
\int_0^{\infty} (1-s+sx^2)x e^{\frac{-x^2}{2}}\erf(x\sqrt{\frac{\kappa}{2}})  +\sqrt{\frac{2}{\pi}}\sigma(s) x^2e^{\frac{-x^2}{2(1-s)}} dx 
\NEQ
\BEQ \label{AM1TIntB}
= (1-s)\frac{\sqrt{\kappa}}{\sqrt{1+ \kappa}}+s\sqrt{\kappa}\frac{2 \kappa +3}{(1+\kappa)^{3/2}} + \sigma(s)*(1-s)^{3/2} = 2 \sqrt{s}  \ .
%e^{\frac{-h^2}{2(1-\theta)}}
%h- \sqrt{\frac{2(1-\theta)}{\pi}}\left[\tan^{-1}(\sqrt{\kappa}) + s \frac{\sqrt{\kappa}}{\kappa +1} \right] + \sqrt{\frac{2s(1-s)^3(1-\theta)}{\pi}} %e^{\frac{-h^2}{2(1-\theta)}}
\NEQ
%\BEQ \label{AM1TIntC}
%= G_{1,2}(s)\equiv (1-s)^2 (\sqrt{s}+\sqrt{\kappa}) + s \sqrt{\kappa}(3-s).
%\NEQ
%where $\kappa = s/(1-s)$ and $\sigma(s)\equiv \sqrt{s(1-s)}$.
We used the identities:
$\int_0^{\infty} x\ \erf(\frac{ax}{\sqrt{2}}) \exp(\frac{-x^2}{2})dx =\frac{a}{\sqrt{1+ a^2}}$ and
$\int_0^{\infty}x^3 \erf(\frac{ax}{\sqrt{2}}) \exp(\frac{-x^2}{2})dx =$ $a(2 a^2+3)/ (a^2 +1)^{3/2}$.

\section{Appendix B: Details of Brownian Simulations \label{AppSimDesc}}

In the simulation, we generate about  $17$ million Brownian paths.
We then bin the paths in a ``grid'' of bins in $(close, \theta, high)$ space.
We use quantile-like coordiantes in the $(q_1,q_2, q_3)$ parameterization.
In the first phase space direction, $q_1=c$,
compute bin boundaries so that the number of curves are roughly equal in each bin.
For each one dimensional bin, compute bin boundaries in the second coordinate direction, $h$, so that the number of bins
is roughly equal. Finally, for each of the two dimensional bins, compute bins in the third direction, $\theta$.
Often, we choose a grid in $c$ to correspond to quantiles of the normal distribution.
(We center the $c$ values to be the midpoints of the $c$ grid.)
The advantage of this approach is that the number of curves in each simulations is approximately the same.
The disadvantage of this approach is that a convergence analysis is hard because the coordinates vary from simulation
to simulation and halving the coordinate distance cannot be done in this adaptive gridding.

Given an ensemble of Brownian paths, $\{B_i(t)\}$, we can create an equivalent ensemble of Brownian paths, $\{B_i(t,c)\}$, with right endpoint $c$, using the formula: $B_i(t,c) \equiv B_i(t) - (B_i(t=1) -c) t$. In our case, we transform the original 10 to 36 million Brownian paths for each value of $c$, calculate the maximum, $h$, and arg-maximum, $\theta$, of $B_i(t,c)$. We then bin the paths into an $nbins \times nbins$ grid in  $(\theta,\ h)$ space based on the quantiles of the $(\theta,\ h)$. A typical values of $nbins$ are from $80$ to $120$.
For each bin in the $(close,\ \theta,\ h)$ grid, we calculate the mean and variance for each time and bin.

To describe these simulations, we call $nSteps$, the number of timesteps in the any one path simulation and $nSim$ is the number of Brownian paths that we compute in the simulation before shifting the paths to have final value $c$.

%To create the bins, we originally quantiles of $\theta$ and $h$. However this results in very wide buckets at the end. We modify this to choose
%buckets so the integral of $density^{.25}$ is constant in each bucket. We also find that
The values of $\argmax$ occur only on the time grid and
this makes the quantile buckets discrete. To regularize the discreteness of $\argmax$, we use a quadratic interpolation of $\argmax$, i.e.
for a given realization of $B(t)$ on a discrete path, we find the $\argmax$ on the grid, $\{t_i \}$. We then use quadratic intepolation of $B(t)$
on $B(t_{i-1})$, $B(t_i)$ and $B(t_{i+1})$ to find a better approximation to the continuous time value of $\argmax$.
A more accurate calculation is to use the expection of the $\argmax$ conditional on the three point values $B(t_{i-1})$, $B(t_{i})$, $B(t_{i+1})$ where the discrete minimum occurs at $B(t_i)$ \cite{AC}.

\section{Appendix C: Plots of Mean and Variance Given High and Argmax \label{EVHA}}

\begin{figure}[htbp]
\centering
\includegraphics[scale=\FSCL]{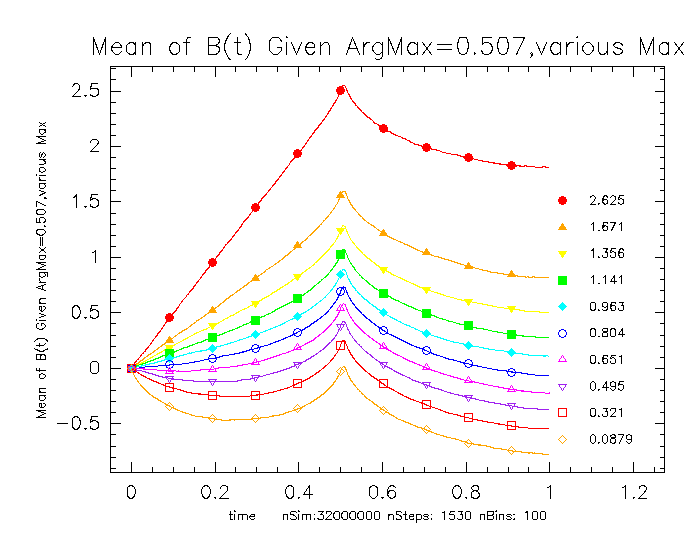}
\caption{$E[B(t |argmax=0.507, various\ high)]$. The $h$ values are spaced at roughly quantile values of $k/9$ given.
 The values of the high are given in the legend.}
\label{fi:AHMeanGiven2ArgMax=0.507AH}
\end{figure}

\begin{figure}[htbp]
\centering
\includegraphics[scale=\FSCL]{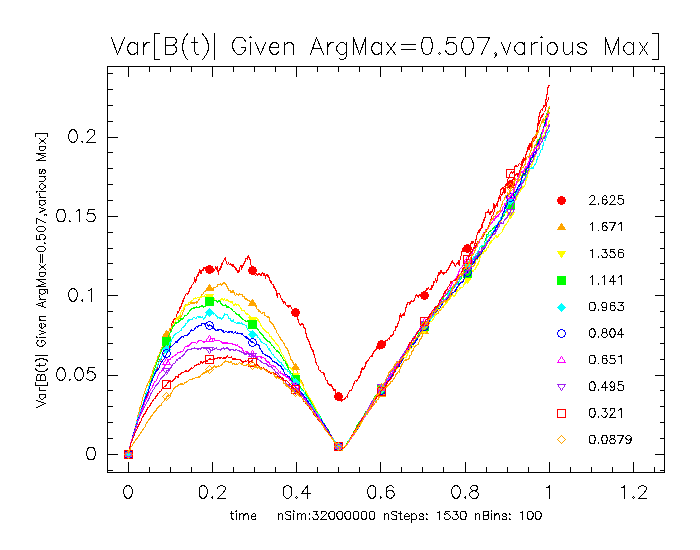}
\caption{$Var[B(t |\argmax=0.507, various\ high)]$ where the values of $\max$ as given in the legend.}
\label{fi:AHVarGiven2ArgMax=0.507AH}
\end{figure}

\begin{figure}[htbp]
\centering
\includegraphics[scale=\FSCL]{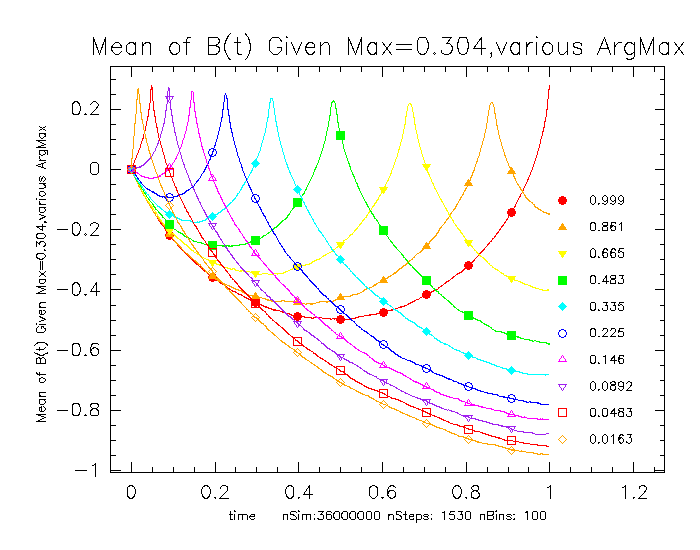}
\caption{$E[B(t |max=0.304, \argmax)]$. Here .304 is roughly the .2 quantile of the high. The $\theta$ values are spaced at roughly quantile values of $k/9$ given, $h$.}
\label{fi:HAMeanGiven2Max=0.304HA}
\end{figure}

\begin{figure}[htbp]
\centering
\includegraphics[scale=\FSCL]{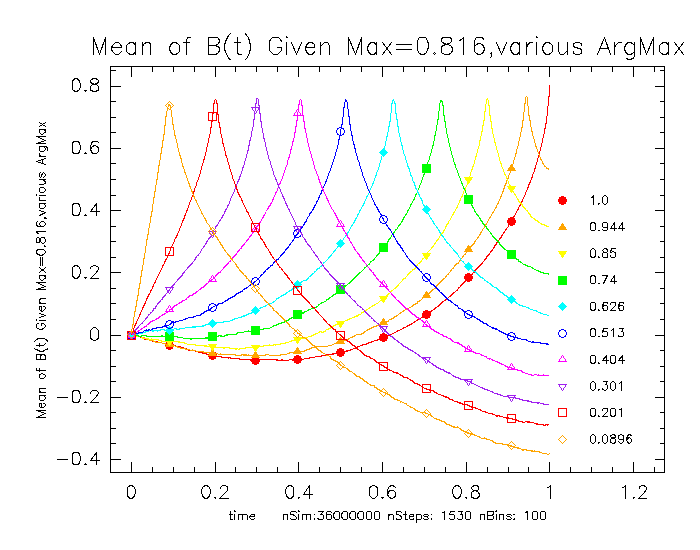}
\caption{$E[B(t |max=0.816, various\ \argmax)]$. Here .816 is roughly the median of the high. The $\theta$ values are spaced at roughly quantile values of $k/9$ given, $h$.}
\label{fi:HAMeanGiven2Max=0.816HA}
\end{figure}

\begin{figure}[htbp]
\centering
\includegraphics[scale=\FSCL]{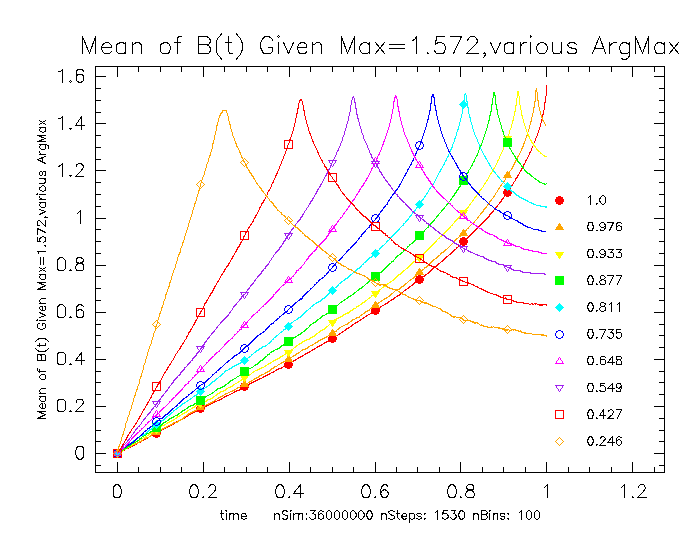}
\caption{$E[B(t |max=1.572, various\ \argmax)]$. Here 1.572 is roughly the .2 quantile of the high. The $\theta$ values are spaced at roughly quantile values of $k/9$ given, $h$.}
\label{fi:HAMeanGiven2Max=1.572HA}
\end{figure}

Figure \ref{fi:AHMeanGiven2ArgMax=0.507AH} displays the expectation of $B(t|\theta=.507, h)$ for various values of the high.
The maximum is at $\theta$ and decays as $t \rightarrow 1$. The values of the high are spaced accordingly:
For the given values of $\theta$, we display the the values of the high corresponding to the second smallest and second largest bins and eight equi-spaced values of $h$.  
Figure \ref{fi:AHMeanGiven2ArgMax=0.507AH} displays the expectation of $B(t|\theta=.507, h)$ for various values of the high and Figure \ref{fi:AHVarGiven2ArgMax=0.507AH} displays the analogous variance.
The variance should be zero at $\theta$ however the curve for the largest value of $h$ does not satisfy this as the discrete binning effect has curves corresponding to a variety of different $\theta$ and the mean of the curves in this edge bin
is not equal to the average of the expectations but simply the average of the realizations.
We see that the variance is roughly parabolic for $t \in [0,\theta]$ and then grows beyond $\theta$.
Figures \ref{fi:HAMeanGiven2Max=0.304HA}-\ref{fi:HAMeanGiven2Max=1.572HA} display the other cross sections for different values of $\theta$ with a fixed maximum, $h$.
Plotting the variance for these cross sections again shows a rounded shape for $t \le \theta$ and a roughly linear growth rate for $t \ge \theta$. For $t<\theta$, the curves rise more rapidly near $t=0$ than they decrease near $t\approx \theta$.
The local maximum of the variance occurs before $t=\theta/2$.
%\begin{figure}[htbp]\centering
%\includegraphics[scale=\FSCL]{{VarGiven2Max=0.304HA_sim36000000_stp1530_Bn100}.png}
%\caption{$Var[B(t |close=VarGiven2Max=0.304HA, high=, low)]$:VarGiven2Max=0.304HA}
%\label{fi:HAVarGiven2Max=0.304HA}
%\end{figure}
\begin{figure}[htbp]
\centering
\includegraphics[scale=\FSCL]{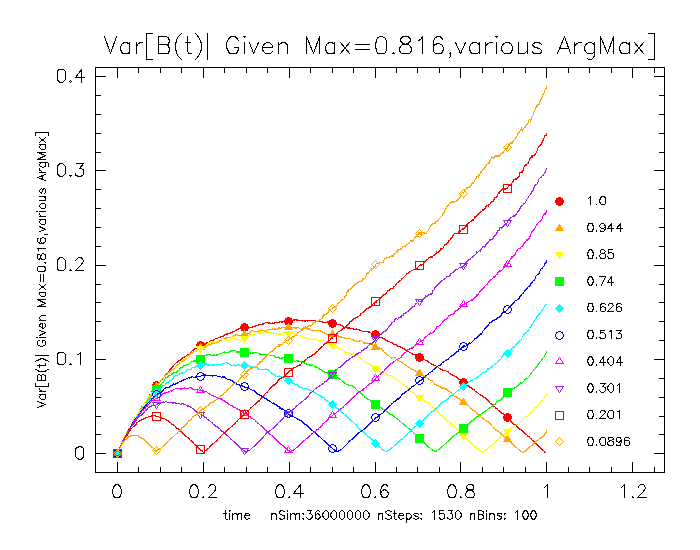}
\caption{$Var[B(t |max=0.816, various\ \argmax)]$}
\label{fi:HAVarGiven2Max=0.816HA}
\end{figure}
%\begin{figure}[htbp] \centering
%\includegraphics[scale=\FSCL]{{VarGiven2Max=1.572HA_sim36000000_stp1530_Bn100}.png}
%\caption{$Var[B(t |close=VarGiven2Max=1.572HA, high=, low)]$:VarGiven2Max=1.572HA}
%\label{fi:HAVarGiven2Max=1.572HA}
%\end{figure}

\begin{figure}[htbp]
\centering
\includegraphics[scale=\FSCL]{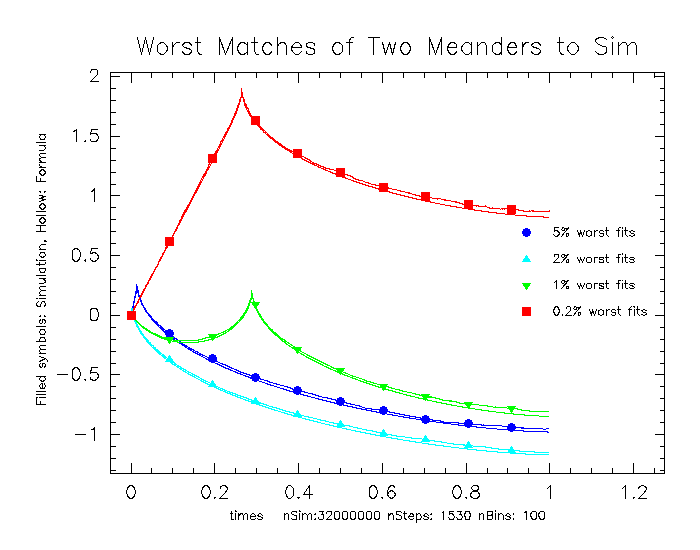}
\caption{Comparison of  Formula and Simulation:\n
Blue: $5\%$ worst MSE Mean:0.000433 high:0.267, argmax:0.0142\n
Cyan: $2\%$ worst MSE Mean:0.000554 high:0.0825, argmax:0.00224 \n
Green: $1\%$ worst MSE Mean:0.000639 high:0.21, argmax:0.289\n 
Red: $0.2\%$ worst MSE Mean:0.000826 high:1.897, argmax:0.265\n
 }
\label{fi:AHMeanCombo}
\end{figure}

\begin{figure}[htbp]
\centering
\includegraphics[scale=\FSCL]{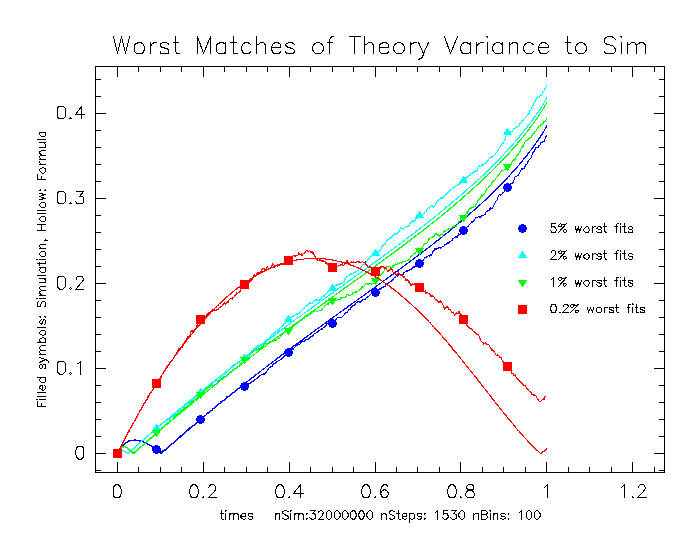}
\caption{Comparison of Meander Variance Formula and Simulation:\n
  Blue: $5\%$ worst MSE Var:.0000541, high:0.329, argmax:0.102\n
Cyan:$2\%$ worst MSE Var:0000943, high:0.127, argmax:0.0247\n
Green:$1\%$ worst MSE Var:0.000143, high:0.46, argmax:0.0376 \n
Red:$0.2\%$ worst MSE Var:0.000857, high:3.658, argmax:0.986\n}
\label{fi:AHVarCombo}
\end{figure}

Figure \ref{fi:AHMeanCombo} and Figure \ref{fi:AHVarCombo} compare our simulations with the analytic expressions in \eqref{EHTh1} - \eqref{EHTh2}
In the simulation, we use $120^3 \approx 1.7M$ bins for binning the values of $(close, \theta, high)$.
For each bin, we compute the expectation and variance using \eqref{MnCThHa} - \eqref{VrCThHb}, and then evaluate the
time average of the mean square error. We then sort the bins from worst to best MSE. Not surprisingly
the worst fitting bins often occur where the bins are the largest and the bias error, from having curves
with slightly different parameters, is the largest. Figure \ref{fi:AHMeanCombo} plots the curves that correspond to
the worst $5\%$, $2\%$, $1\%$ and $0.2\%$. Since we evaluate the worst fitting curves for the expectation and
variance separately, the parameters for the worst fits for the expectation differ from the corresponding parameters for the vairance plots. The curves overstruck by symbols are the simulation curves. The analytic formula curves have the same color but no symbol.
Overall, the concatenated meander expectation fits the simulation expectation very well.
%To further illustrate this, we overlay all four sets of curves on the same plot.

The MSE in Figure \ref{fi:AHMeanCombo} is less than the MSE in Figure \ref{fi:MeanderMeanCombo}.
%we also dispalay the analogous plots for $E[B(t|\theta, h)$.
As expected, averaging over values of $c$ reduces the noise in the simulations and reduces the discrepancy.

We conclude with a comparison of the variance of the simulation with the concatenated meander in \eqref{VrCThHa} - \eqref{VrCThHb}.

\section{Appendix D: Time Dependent Behavior of Moments Given Close, Argmax, High \label{EVCAH}} 

The figures in this section are from our simulation without the corresponding
analytic calculations.

\subsection{Expectation Given Close, Argmax, High \label{ECAH}} 
The next five figures display $E[B(t |close, \argmax, high)]$. Each figure fixes values of the close and $\theta= \argmax\{B(t)\}$.
The values of the close are roughly $\{-1,0,+1\}$. In this section, we use 12.5 million realizations with 1850 time steps.
We divide the parameter space into 120 bins in parameter direction. 
For each realization, we shift all of the sample paths to have the same close. For each close, the 12.5M realizations
are binned into the $120^2$ bins. Thus the empirical mean and variance are noisier than in Section \ref{EVCA} where the realizations
were divided into only $120$ bins. We plot both the empirical cuves from our simulation and the theoretical values from
\eqref{MnCThHa} - \eqref{VrCThHb} in the same color. The curves overstruck by symbols are the simulation curves. The analytic formula curves have the same color but no symbol. 
The smooth curves are the theoretical ones while the noisy curves are the simulated ones.
\begin{figure}[!htbp]
\centering
\includegraphics[scale=\FSCL]{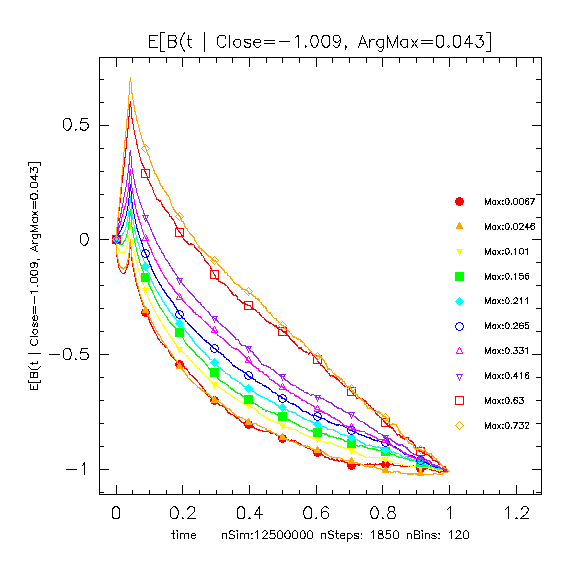}
\caption{$E[B(t |close, \argmax, high)]$ where close$=-1.009$ and argmax$=0.043$.
  Each curve is a given value of high$=max\{B(t)\}$.}
\label{fi:CAHMeanGivenClose:-1.009 ArgMax:0.043CAH}
\end{figure}
\begin{figure}[htbp]
\centering
\includegraphics[scale=\FSCL]{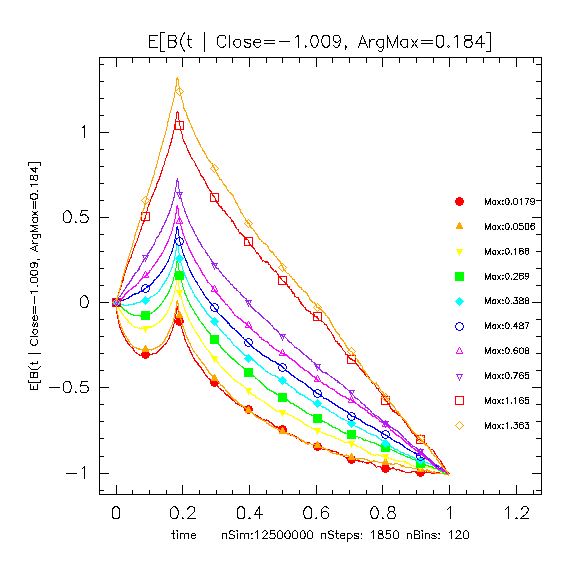}
\caption{$E[B(t |close, \argmax, high)]$ where close$=-1.009$ and $\theta=0.184$. The values of the high are given in the legend.
  For larger values of the high, the $E[B(t)]$ is approximately piecewise linear with a maximum at $\theta$.}
\label{fi:CAHMeanGivenClose:-1.009ArgMax:0.184}
\end{figure}
\begin{figure}[htbp]
\centering
\includegraphics[scale=\FSCL]{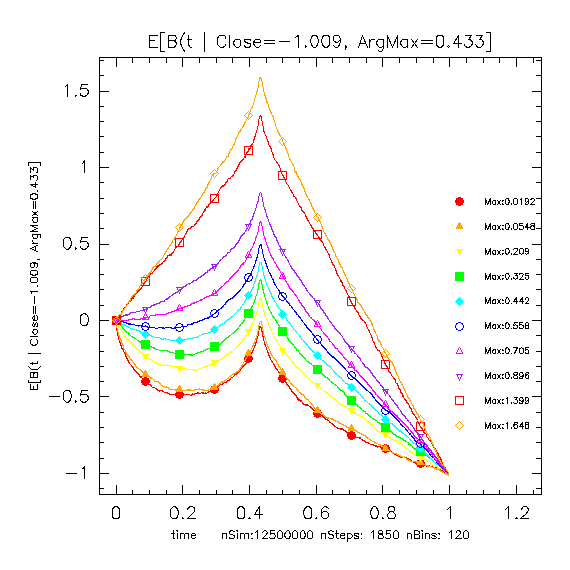}
\label{fi:CAHMeanGivenClose:-1.009ArgMax:0.433}
\caption{$E[B(t |close, \argmax, high)]$ where close$=-1.009$ and argMax$=0.433$.
  Each curve is a given value of high$=max\{B(t)\}$ as denoted in the legend.}
\end{figure}

\begin{figure}[htbp]
\centering
\includegraphics[scale=\FSCL]{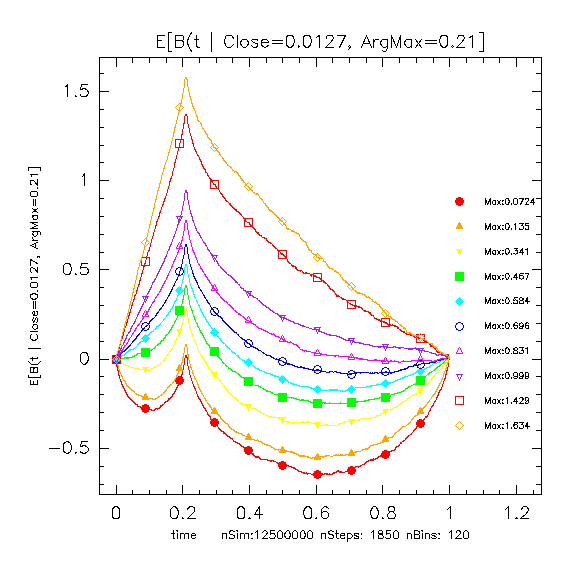}
\label{fi:CAHMeanGivenClose:0.0127 ArgMax:0.21CAH}
\caption{$E[B(t |close, high, \argmax)]$ where close$=0.0127$ and $\theta=0.21$.
 For larger values of the high, the $E[B(t)]$ grows approximately linearly until $\theta$ and then decays back to its final value.}
\end{figure}
For $close=0$, we have the symmetry: $E[B(t |c=0, h, \theta)]=E[B(t |c=0, h,1- \theta)] $.
\begin{figure}[htbp]
\centering
\includegraphics[scale=\FSCL]{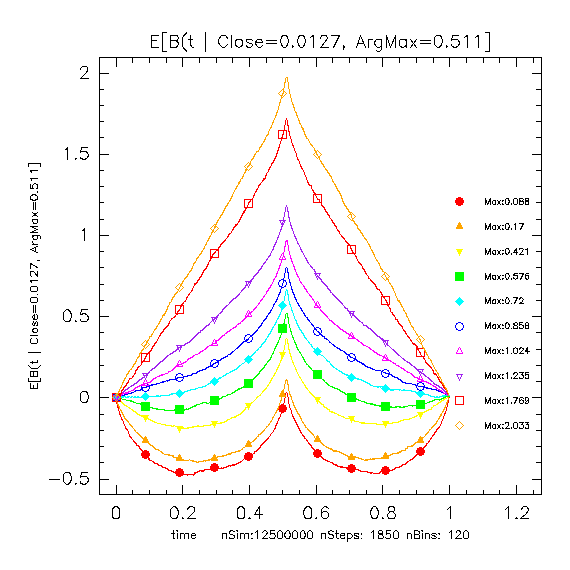}
\label{fi:CAHMeanGivenClose:0.0127 ArgMax:0.511}
\caption{$E[B(t |close, \argmax, high)]$ where close$=0.0127$, $\theta=0.511$, and $\max {B}$ is given in the legend. .
Since the close is nearly zero and $\theta$ is nearly $0.5$, the curves are nearly symmetric in time.}    
  %Each curve is a given value of high$=max\{B(t)\}$. 
\end{figure}
Rather than plot figures for $close=1$, we again the reflection symmetry:
$E[B(t |-c, \theta, h)] = E[B(1-t |c,1- \theta,h+c)] - c$.

\clearpage
\subsection{Variance Given Close, Argmax, High \label{VCAH}} 
We now display $Var[B(t |close, \argmax, high)]$. Each figure fixes values of the close and $\theta= \argmax\{B(t)\}$. The curves overstruck by symbols are the simulation curves. The analytic formula curves have the same color but no symbol. The simulation curves have a larger value because they include the $bias^2$ term.
\begin{figure}[!htbp]
\centering
\includegraphics[scale=\FSCL]{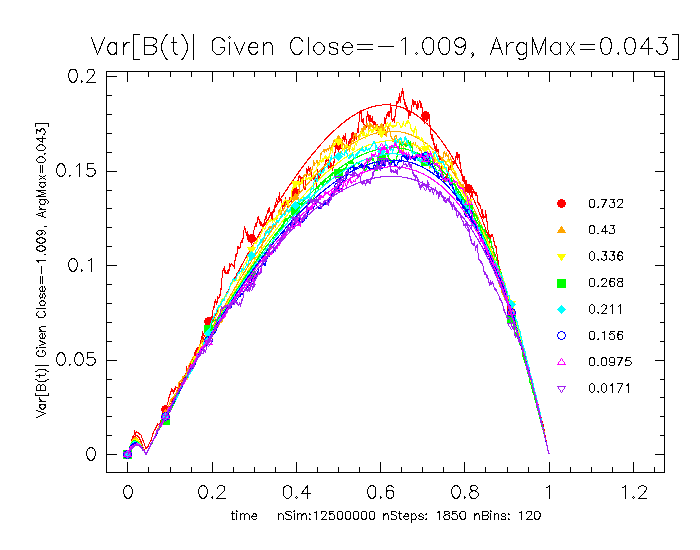}
\caption{$Var[B(t |close, \argmax, high)]$ where close$=-1.009$ and argMax$=0.043$. The value of $\theta$ is small,
  therefore the variance has a local minimum at $\theta$ and then increases. The variance profile appears especially  flat near its maximum.}
%  Each curve is a given value of high$=max\{B(t)\}$ }
\label{fi:VarGivenClose:-1.009ArgMax:0.043CAH}
\end{figure}
\begin{figure}[htbp]
\centering
\includegraphics[scale=\FSCL]{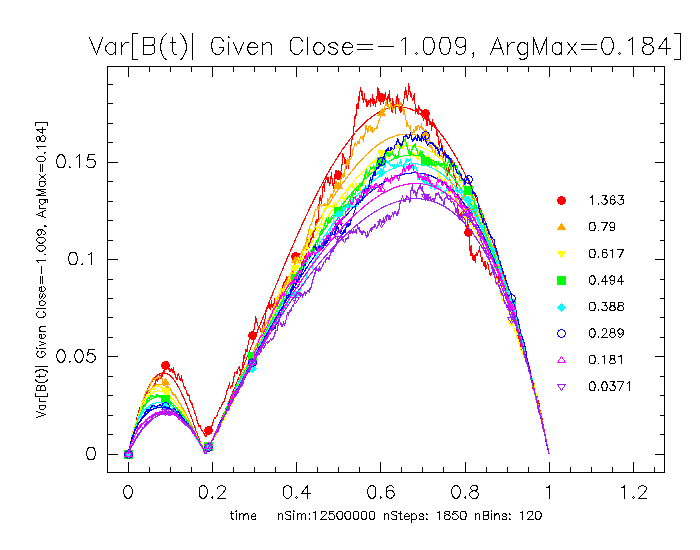}
\caption{$Var[B(t |close, \argmax, high)]$ where close$=-1.009$ and argMax$=0.184$.
  Each curve is a given value of high$=max\{ B(t)\}$. The curves are bimodal with a local minimum at $t=\theta$.}
\label{fi:VarGivenClose:-1.009ArgMax:0.184}\
\end{figure}
\begin{figure}[htbp]
\centering
\includegraphics[scale=\FSCL]{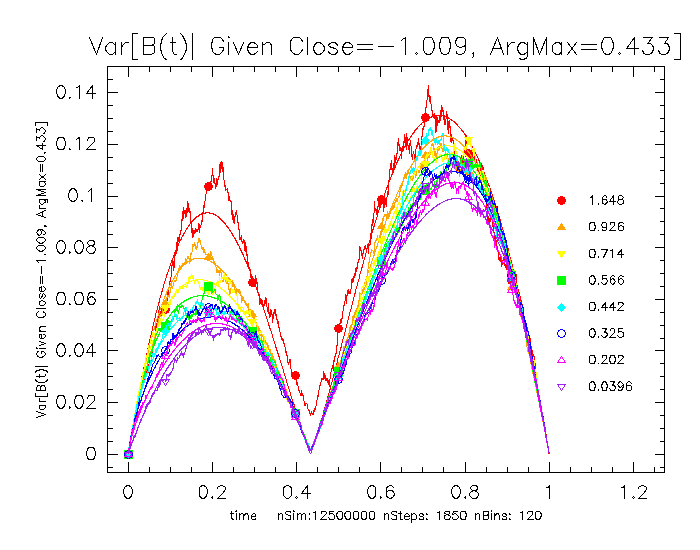}
\caption{$Var[B(t |close, \argmax, high)]$ where close$=-1.009$ and argMax$=0.433$.
  The maximum of the variance is smaller when $\theta$ is larger. Since the variance is small for $t$ near $\theta$, the variance cannot grow as large when $\theta$ splits time into two roughly equal time intervals.}
\label{fi:VarGivenClose:-1.009ArgMax:0.433CAH}
\end{figure}
Rather than plot figures for $close=1$, we use the reflection symmetry:
$Var[B(t |-c, \theta, h)] = Var[B(1-t |c,1- \theta,h+c)]$.
%Figure \ref{fi:CAH:VarGivenClose=1.046} and figure \ref{fi:CAH:VarGivenClose=-1.009} have reflection symmetry:
\begin{figure}[htbp]
\centering
\includegraphics[scale=\FSCL]{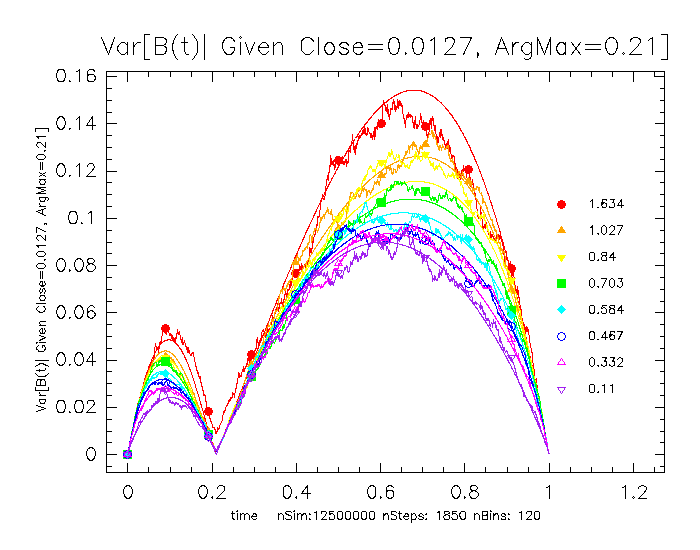}
\caption{$Var[B(t |close=0.0127, \theta=0.21, high)]$ for various $h$ values.
  Each curve is a given value of high$=max B(t)$ as denoted in the legend.}
\label{fi:VarGivenClose:0.0127ArgMax:0.21CAH}
\end{figure}
For $close=0$, we have the symmetry: $Var[B(t |c=0, \theta,h)]=Var[B(t |c=0, 1- \theta,h)] $
\begin{figure}[htbp]
\centering
\includegraphics[scale=\FSCL]{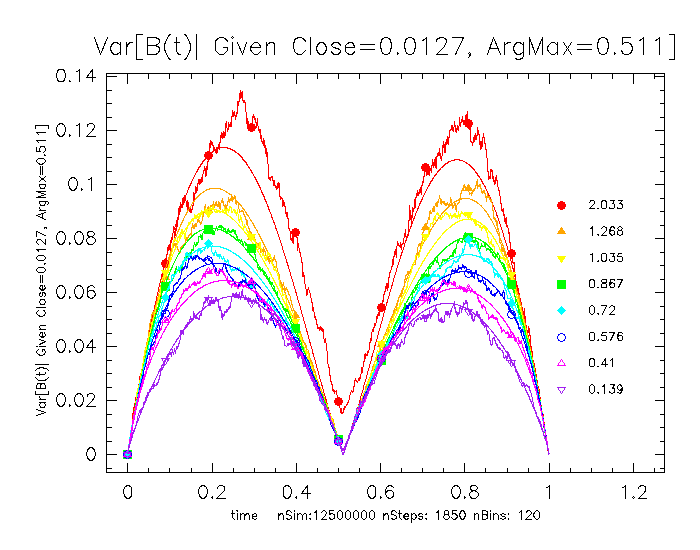}
\caption{$Var[B(t |close=0.0127,\theta=0.511,high)]$. The variance is nearly symmetric in time.}
%    Each curve is a given value of high$=max B(t)$. }
\label{fi:VarGivenClose:0.0127ArgMax:0.511}
\end{figure}

\section{Appendix E: Mean and Variance Given Final Value and Argmax \label{EVCA}}
%\afterpage{\FloatBarrier}

In the last section, we integrated over the $\close$ to calculate the moments conditional on $(high, \argmax)$.
We now display the mean and variance of $B(t|\theta,c)$.
Similar but more complicated moment calculations are possible for the case of the previous section where $(\theta,c)$ are given and one integrates over $h$. The figures in this section are from our simulation without the corresponding
analytic calculation.

\subsection{Time Dependence of Expectation given Close and Argmax} 

\begin{figure}[htbp]
\centering
\includegraphics[scale=\FSCL]{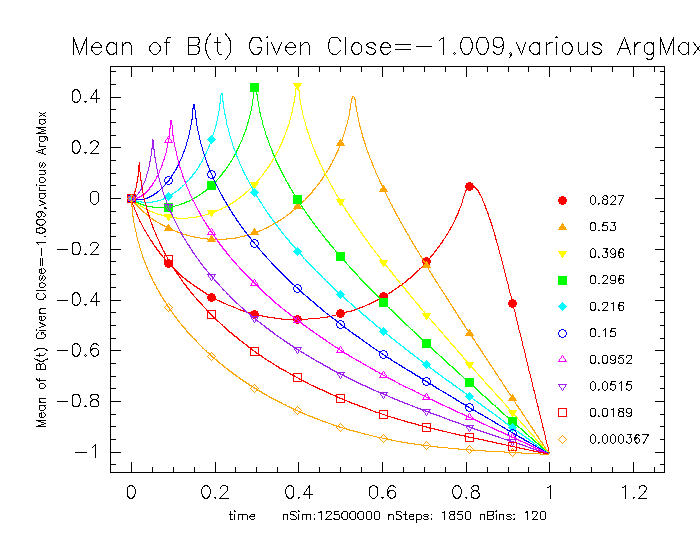}
\caption{$E[B(t |close, \argmax)]$ where close$=-1.009$. Each curve is a given value of $\theta=\argmax\{B\}$.
  Since the close is roughly at $- \sigma$, the argmax is preferentially located near zero. The expectation is peaked at $t=\theta$.
The values of $theta$ are given on the legend.}
%The symbol markers are placed at equal distiances in bin number.
\label{fi:CAH:MeanGivenClose=-1.009}
\end{figure}

\begin{figure}[htbp]
\centering
\includegraphics[scale=\FSCL]{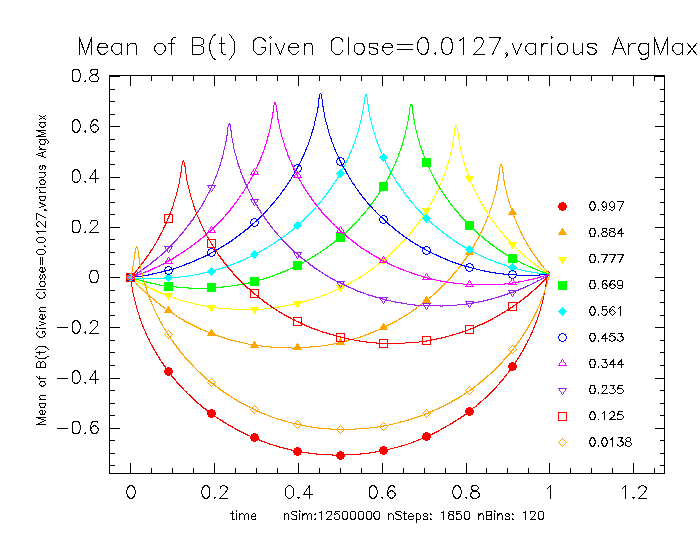}
\caption{$E[B(t |close, \argmax)]$ where close$=0.0127$. Each curve is a given value of $\theta=\argmax\{B\}$.
Since the close is near zero, if $\theta$ is small, the expectation after the zero-crossing resembles a semicircle.}
%The symbol markers are placed at equal distiances in bin number.}
\label{fi:CAH:MeanGivenClose=0.0127}
\end{figure}

\begin{figure}[htbp]
\centering
\includegraphics[scale=\FSCL]{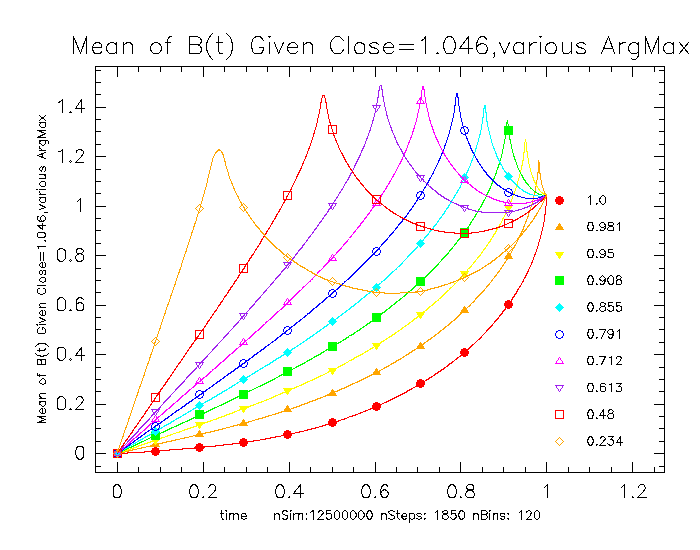}
\caption{$E[B(t |close=1.046, \argmax)]$. Each curve is a given value of $\theta=\argmax\{B\}$ as given on the legend.}
%The symbol markers are placed at equal distances in bin number.}
\label{fi:CAH:MeanGivenClose=1.046}
\end{figure}

Figures \ref{fi:CAH:MeanGivenClose=-1.009}, \ref{fi:CAH:MeanGivenClose=0.0127} and \ref{fi:CAH:MeanGivenClose=1.046}
display the expectation of $E[B(t |close, \argmax)]$. Each figure fixes a value of the close.
The values of the close are roughly $\{-1,0,+1\}$. The simulations use 12.5 million realizations with 1850 time steps.
For each realization, we shift all of the sample paths to have the same close. For each close, the 12.5M realizations
are binned into the 120 bins.

Figure \ref{fi:CAH:MeanGivenClose=1.046} and Figure \ref{fi:CAH:MeanGivenClose=-1.009}  have reflection symmetry:
$E[B(t |-c, \theta)] = E[B(1-t |c,1- \theta)] - c$. This is a consequence of the reversability of Brownian motion.

%\clearpage
\subsection{Time Dependence of Variance given Close and Argmax} 
\begin{figure}[htbp]
\centering
\includegraphics[scale=\FSCL]{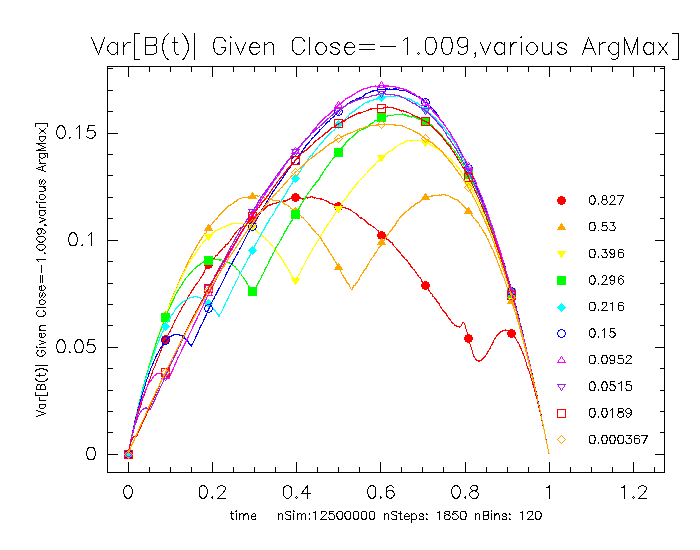}
\caption{$Var[B(t |close=-1.009, \argmax)$ where close$=-1.009$ whith the values of $\theta$ given on the legend.}
\label{fi:CAHVarGiven2Close=-1.009}
\end{figure}

\begin{figure}[htbp]
\centering
\includegraphics[scale=\FSCL]{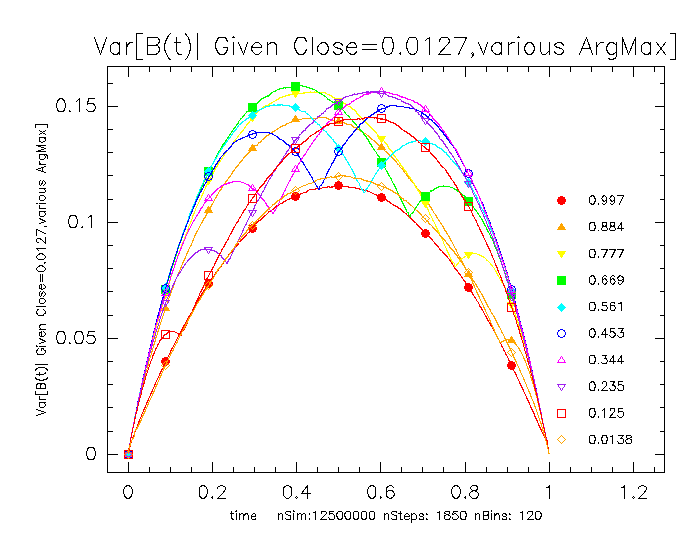}
\caption{$Var[B(t |close=.0127, \argmax)$ where close$=0.0127$.}
\label{fi:CAHVarGiven2Close=0.0127}
\end{figure}

\begin{figure}[htbp]
\centering
\includegraphics[scale=\FSCL]{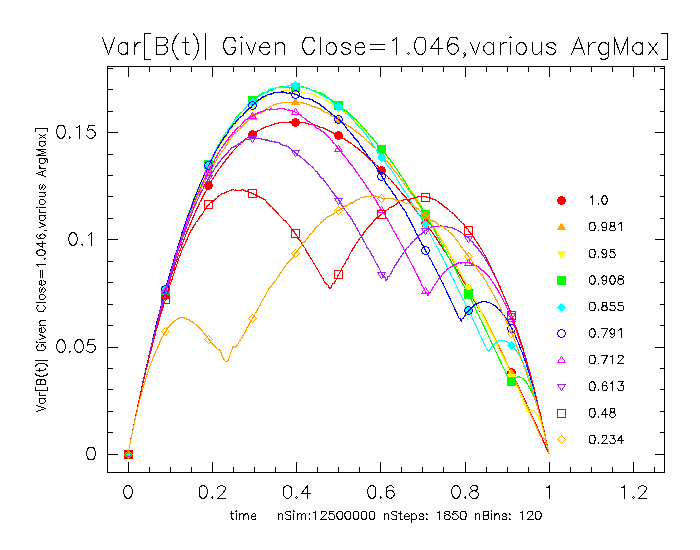}
\caption{$Var[B(t |close=1.046, \argmax)$ where close$=1.046$.}
\label{fi:CAHVarGiven2Close=1.046}
\end{figure}

Figures \ref{fi:CAHVarGiven2Close=-1.009}, \ref{fi:CAHVarGiven2Close=0.0127} and \ref{fi:CAHVarGiven2Close=1.046}
display the variance of $B(t |close, \argmax)$. Each figure fixes a value of the close.
The values of the close are roughly $\{-1,0,+1\}$.

Figure \ref{fi:CAHVarGiven2Close=1.046} and Figure \ref{fi:CAHVarGiven2Close=-1.009} have reflection symmetry:
$Var[B(t |-c, \theta)] = Var[B(1-t |c,1- \theta)]$.

\ACKNO{The author thanks his referees for their helpful remarks.}
%\ACKNO{On behalf of all authors, the corresponding author states that there is no conflict of interest.}

\end{document}